\documentclass[12pt]{article}
\usepackage{amsmath,amsthm,verbatim,amssymb,amsfonts,amscd,diagrams, graphics}
\theoremstyle{plain}
\newtheorem{theorem}{Theorem}[section]
\newtheorem{corollary}[theorem]{Corollary}
\newtheorem{lemma}[theorem]{Lemma}
\newtheorem{proposition}[theorem]{Proposition}

\theoremstyle{definition}
\newtheorem{definition}[theorem]{Definition}

\theoremstyle{remark}
\newtheorem{remark}[theorem]{Remark}

\newarrow{ul}---->
\newarrow{Backwards}<----

\newcommand{\td}[1]{\tilde{#1}}
\newcommand{\into}{\hookrightarrow}
\newcommand{\Z}{\mathbb{Z}}

\newcommand{\R}{\mathbb{R}}

\newcommand{\N}{\mathbb{N}}
\newcommand{\bd}{\partial}

\renewcommand{\H}{\mathbb H}

\newcommand{\mc}[1]{\mathcal{#1}}

\newcommand{\holink}{\text{holink}}

\newcommand{\hl}{\operatorname{holink}}

\newcommand{\mf}{\mathfrak}

\newarrow{Onto}----{>>}
\newarrow{Equals}=====

\begin{document}

\title{Intersection homology of stratified fibrations and neighborhoods}
\author{Greg Friedman\\Texas Christian University\\Dept. of Mathematics\\\\g.friedman@tcu.edu }
\date{November 18, 2005}
\maketitle

\begin{abstract}
We derive spectral sequences for the intersection homology of stratified fibrations and approximate tubular neighborhoods in manifold stratified spaces. These neighborhoods include regular neighborhoods in PL stratified spaces. 
\end{abstract}


\tableofcontents

\section{Introduction}
As with other homology theories, to compute intersection homology one first needs to be able to compute on the most fundamental spaces. One might then proceed from local to global computations by applying inductively such machinery as long exact sequences or spectral sequences. While it may only rarely be possible to obtain complete information by these methods, one hopes at least to be able to extract some useful invariant information, such as the existence or limitations on certain types of torsion or constraints on Betti numbers. Such programs are carried out in special cases, for example, in \cite{GBF2} and \cite{LM06} to obtain properties on  intersection homology  invariants of knots and hypersurfaces. 

In the theory of stratified spaces, the most local information lives in \emph{distinguished neighborhoods} of points - neighborhoods of  (or stratum-preserving homotopy equivalent to) the form $\R^k\times cL$, where $L$ is a lower-dimensional stratified space and $cL$ is the open cone on it. So the first computations in the theory focused on formulae for intersection homology of cones, suspensions, and products, the last usually with respect to euclidean spaces or manifolds. And indeed, on sufficiently well-behaved spaces, such as topological pseudomanifolds, the intersection homology modules of these local structures governs the global modules via sheaf theoretic spectral sequence $H^*(X;\mc I^{\bar p} \mc H^*)\Rightarrow I^{\bar p}H^*(X)$. Here the first term is the cohomology of $X$ with sheaf coefficients in $\mc{IH}$, whose stalks are the intersection cohomology of the  links $L$. The second term is the intersection cohomology of $X$. In fact, the axiomatic characterization of the intersection chain sheaf on pseudomanifolds, proven in \cite{GM2} (and presented in a nice exposition in \cite{Bo}), demonstrates that the cone and product formulas completely govern the behavior of intersection homology. This can also be seen in King's proof of the topological invariance of singular chain intersection homology \cite{Ki}, in which most arguments ultimately come down to the behavior of intersection homology on cones, products (the K\"unneth theorem), and unions (the Mayer-Vietoris sequence). 

Of course stratified spaces possess natural intermediate structures between the local and global, namely the strata themselves or, more generally, closed unions of connected components of strata, called \emph{pure subsets}. Unfortunately, these subspaces cannot themselves be mined for intersection homology information, as the defining properties of intersection chains (or intersection chain sheaves) is that they live primarily on the dense open top stratum - only dipping into the lower strata enough to keep things interesting. So to study this intermediate structure intersection-homologically, we must instead study neighborhoods, most naturally regular neighborhoods of pure subsets and their analogues in broader topological and homotopical categories (see Section \ref{S: neighborhoods}, below, for precise details on the types of neighborhoods we consider). Of course, the ordinary homology of such neighborhoods is just the homology of the ``base'' (the space whose neighborhood we study), because the neighborhoods and their bases are homotopy equivalent. But intersection homology is not a homotopy invariant - it is only invariant under stratum-preserving homotopy equivalences (and then only with compact supports)- so we expect the intersection homology of the neighborhoods to contain homological information reflecting both the homological character of the base space and of how it fits into the stratified structure of the space as a whole. 

It has long been natural to study neighborhoods of subspaces in conjunction with bundle theories, perhaps growing from the early identification of tubular neighborhoods of smooth submanifolds with normal vector bundles. Similar correspondences then worked their way into more general categories with block bundles, microbundles, etc. Even when a purely geometric bundle structure is not available, there is hope that one can employ homotopy methods to identify neighborhoods with fibrations, and then one might hope to attack homological questions with generalizations of the Leray-Serre spectral sequence. In particular, the idea of identifying neighborhoods with path space fibrations goes back at least to Fadell \cite{Fad65}, and it is built into the definition of Quinn's manifold homotopically stratified spaces (MHSSs) \cite{Q1}.

The most natural neighborhoods of strata (or pure subsets) in MHSSs turn out, in fact, not to have natural bundle structures, but they are stratum-preserving homotopy equivalent to stratified fibrations. These stratified fibrations possess several of the  properties one would expect of something called a fibration (see \cite{Hug,GBF3}), and, in particular, there is a stratified local homotopy uniformity over strata of the base space. Thus it is reasonable to expect that one might exhibit a spectral sequence that abuts to the intersection homology of the neighborhood and whose $E^2$ terms are written in terms of the (co)homology of the base and with coefficients that are locally-constant on each stratum and that reflect the local intersection homology  at points of those strata.

Since we thus will have coefficient systems varying from stratum to stratum on the base, the most natural language in which to organize this data is that of sheaf cohomology. In these terms, we want to find a complex of constructible sheaves on the base, cohomologically locally constant on each stratum of the base, and  the cohomology modules of these sheaves should be the  $E^2$ terms of a spectral sequence that abuts to the intersection homology of the neighborhood. 

At this point, we must consider the technical issues. The most obvious thing to try would be to start with the intersection chain sheaf on the neighborhood, as defined in \cite{GM2} for pseudomanifolds and \cite{GBF10} for manifold homotopically stratified spaces (MHSSs), and to use a pushforward to project it down to the base. Then, letting $N$ stand for the neighborhood and $Y$ for the base space, we certainly would have $IH^*(N)=\H^*(N;\mc{IC}^*)\cong \H^*(Y; p_*(\mc{IC}^*))$, where $\H$ denotes sheaf hypercohomology. But $N$, itself, is not necessarily a fibration over $Y$ (stratified or not), so we don't expect $p_*(\mc{IC}^*)$ to have the uniformity properties that would make the pushforward stratified cohomologically locally constant, a property which is essential in performing meaningful calculations. 

On the other hand, we will see that we do have a stratum-preserving homotopy equivalence from $N$ to a stratified fibration $p: E\to Y$, and compactly supported intersection homology is a stratified homotopy invariant. But, unfortunately even if we push or pull the intersection chain sheaf from $N$ to $E$, there is no direct correspondence between compact supports in $E$ and those in $N$, and we cannot ignore supports, since noncompactly supported intersection homology is not a stratified homotopy invariant, except for proper homotopies. But we can not expect $N$ to be properly homotopy equivalent to $E$, which is the mapping cylinder of the evaluation of a path space. In fact, $E$ will not be even locally compact, which can create other technical glitches in the most standard sheaf machinery.

So we take instead the following approach. Rather than using the intersection chain sheaf on $N$ or $E$, we construct a sheaf on $Y$ directly, but modeled upon what a Leray-Serre sheaf should look like if $E$ had the properties we wanted. More precisely, we build a sheaf on $Y$ whose sections are intersection chains in $E$, but with certain support conditions necessary to ensure that the cohomology of this sheaf (also with the appropriate supports) will give the compactly supported intersection homology of $E$, which is the same as that of $N$. But since our chains live in $E$, we will be able to use the fibration properties of $E$ to ensure constructibility of the sheaf on $Y$. In this way, we obtain the desired spectral sequences for the compactly supported intersection homology modules $IH^c_*(N)$, $IH^c_*(N-Y)$, and $IH^c_*(N, N-Y)$. 

In fact, the  methods we will investigate also yield a spectral sequence for the closed support intersection homology $IH^{\infty}_*(N)$, which is the version of singular intersection homology in agreement with the Deligne-sheaf  intersection homology of \cite{GM2}. This will follow once we demonstrate that, for appropriate $N$, $IH^{\infty}_*(N)\cong IH^{\infty}_*(N,N-Y)$. This is still not quite the same as $IH^c_*(N,N-Y)$, but we will see that the local cocompactness of $N-Y$ in $N$ is sufficient to show that this is nonetheless isomorphic to $IH^{\infty}_*(E,E-Y)$, which is computable by similar methods to those used to compute $IH^{c}_*(E,E-Y)\cong IH^c_*(N,N-Y)$. 

In forthcoming work, these computations will be utilized as part of a program to relate singular chain intersection homology to the Deligne sheaf construction \cite{GM2, Bo} on manifold homotopically stratified spaces. 

\subsection{Outline}

Let us now outline in more detail the sections of this paper and the results contained therein. 

Section \ref{S: bg} contains a review of the necessary background material, including singular intersection homology, stratified spaces and fibrations, manifold homotopically stratified spaces (MHSSs), homotopy links,  teardrops, and the teardrop topology. 

In Section \ref{S: strat fib}, we study the compactly supported intersection homology of stratified fibrations $p:E\to Y$. Stratified fibrations are important in their own right, and these computations will also be essential for the following study of the intersection homology of neighborhoods. We construct a sheaf $\mc{IZ^*}$ on $Y$ whose hypercohomology with compact supports is $IH^c_*(E)$, and it follows in Theorem \ref{T: strat fib SS} that there is a spectral sequence abutting to $IH^c_*(E)$ whose $E^2$ terms are the sheaf cohomology groups of $Y$ with coefficients in the derived sheaf $\mc{H(IZ^*)}$. We show in Section \ref{S: stalks} how the stalks of the derived sheaf relate to the local fiber intersection homology groups. In particularly, we show that, if $Y$ is a MHSS, the derived sheaf is locally constant on each stratum of $Y$. Throughout, we also consider the  closed support hypercohomology of $\mc{IZ}^*$, which can be interpreted as the ``fiberwise compact'' intersection homology of $E$. Similar results hold, and this version of intersection homology is also needed in our study of the intersection homology of neighborhoods.

Section \ref{S: relative} contains the analogous results concerning relative intersection homology for stratified fibration pairs $p:(E,A)\to Y$.

In Section \ref{S: strat fib closed}, we take up the study of the closed support intersection homology $IH^{\infty}_*(E)$ for a stratified fibration $p:E\to Y$. In Theorem \ref{T: lf IH SS}, we note that there is a spectral sequence for these groups determined by the pushforward $p_*$ of the intersection chain sheaf on $E$. However, we must observe immediately that, without any further restrictions on the spaces, it is difficult to say anything about the cohomology of the stalks of the pushforward because the homotopy properties of the fibration are not immediately compatible with closed support homology theories. Nevertheless, we show in Proposition \ref{P: cocompact} that if $p:(E,A)\to Y$ is a fibration pair such that $p^{-1}(K)\cap (E-A)$ is compact for each compact $K\subset Y$, then $IH^{\infty}_*(E,A)\cong IH^{fc}_*(E,A)$, the fiberwise compact intersection homology studied in Section \ref{S: strat fib}.

Section \ref{S: man base} concerns the special case in which the base of a stratified fibration is an unfiltered manifold. In this situation, it is possible to calculate the stalks of the derived sheaf $\mc{H^*(IZ^*)}$ much more explicitly, and the $E^2$ terms of the previous sections reduce to ordinary homology or cohomology with coefficients in a single local coefficient system. Thus simpler versions of the statements of the previous sections occur and are recorded here. We also show that if stratified fibrations are actually stratified bundles, then similar computations can be made for closed support intersection homology.

In Section \ref{S: neighborhoods}, we begin to apply our results on stratified fibrations to computing the intersection homology of neighborhoods $N$ of pure subsets $Y$ of MHSSs. We first reviews concepts concerning the  nearly stratum-preserving deformation retract neighborhoods (NSDRNs) of \cite{GBF5} and the approximate tubular neighborhoods of Hughes \cite{Hu02}. We show in Corollary \ref{C: NSDRN IHC} that there are spectral sequences for the  compactly supported intersection homology of such neighborhoods $N$, as well as of $N-Y$ and $(N,N-Y)$. The $E^2$ terms are the cohomology modules of $Y$ with sheaf coefficients determined by a stratified fibration over $Y$ coming from the evaluation of stratified path spaces (homotopy links) in $N$.

In Section \ref{S: N closed}, we show that if the NSDRN $N$ satisfies a special condition, which does hold for approximate tubular neighborhoods, then $IH^{\infty}_*(N)$ is homeomorphic to the fiberwise compact intersection homology of a stratified fibration, which implies that there are also spectral sequences, closely related to those described in the preceding paragraph, for computing $IH^{\infty}_*(N)$.

Section \ref{S: pms} concerns pseudomanifolds, the spaces on which intersection homology was originally defined by Goresky and MacPherson \cite{GM1, GM2}. We show that in the context of these spaces, the stalk cohomology of the sheaves arising in the spectral sequences for the computation of intersection homology of neighborhoods can be computed in terms of the geometric links of the strata. Thus, once again, a more concrete situation yields much more specific computational tools

In Section \ref{S: applications}, we outline some elementary applications, and finally, in Section \ref{S: ATN} we prove some of the technical results about approximate tubular neighborhoods that are used earlier in the paper. In particular, we show that approximate tubular neighborhoods are NSDRNs and that they are \emph{outwardly stratified tame}. See Section \ref{S: neighborhoods} for the appropriate definitions.

\begin{remark}[A note on previous work.] The author previously studied the intersection homology of neighborhoods in \cite{GBF3} and \cite{GBF5}; in the interest of clarity, we point out the principle ways in which the current work generalizes that of the earlier papers. 

In \cite{GBF3, GBF5}, an effort was made to avoid sheaf theory in the hopes of obtaining results solely in terms of the ordinary homology of the base space. This required imposing several strong conditions on both the base space and the neighborhoods. In \cite{GBF3}, though manifold homotopically stratified space (MHSS) techniques were employed, the focus was exclusively on regular neighborhoods in stratified PL pseudomanifolds, and the base space was always assumed to be an unfiltered manifold stratum. We obtained spectral sequences whose $E^22$ groups are the ordinary homology of the base with coefficients in a single local coefficient system with fiber intersection homology stalks. 

In \cite{GBF5}, we generalized to neighborhoods of pure subsets in manifold homotopically stratified spaces. However, we required the extra conditions that the base be triangulable and that neighborhoods be ``cylindrical''. The resulting $E^22$ terms are given by the homology of the base in a ``stratified system of coefficients'' whose stalks are locally constant on a given stratum but may vary from stratum to stratum. See \cite{GBF5} for more details. 

In addition, only compactly supported intersection homology is considered in \cite{GBF3, GBF5} 

In this paper, we remove the triangulability condition on the base spaces, which fits more naturally with the category of manifold homotopically stratified spaces, and we remove the cylindricality restriction on neighborhoods. This allows us to consider neighborhoods including Hughes's approximate tubular neighborhoods \cite{Hu02}, the most natural geometric generalization of regular neighborhoods to MHSSs. We also study closed support intersection homology. The cost of the additional generality is the need to use sheaf theoretic techniques and to obtain results in sheaf theoretic formulations. Applications still abound, even in this slightly rarefied setting, but the reader in mind of very concrete computations on triangulated spaces is encouraged to consult \cite{GBF3}, \cite{GBF5}, and \cite{GBF2},  as well. 
\end{remark}

\section{Background and Basic Terminology}\label{S: bg}

\subsection{Intersection homology}

In this section, we provide a quick review of the definition of intersection homology. For more details, the reader is urged to consult King \cite{Ki} and the author \cite{GBF10} for singular intersection homology and the original papers of Goresky and MacPherson \cite{GM1,GM2} and the book of Borel \cite{Bo} for the simplicial and sheaf definitions. Although we do make use of sheaf theory in this paper, there is little direct use of the Goresky-MacPherson-Deligne axiomatic sheaf theoretic formulation of intersection homology.
We shall work with the singular chain intersection homology theory introduced in \cite{Ki} with finite chains (compact supports) and generalized in \cite{GBF10} to include locally-finite but infinite chains (closed supports). 

We recall that this intersection homology  is defined on any filtered space

$$X=X^n\supset X^{n-1}\supset \cdots \supset X^0\supset X^{-1}=\emptyset.$$
In general, the superscript ``dimensions'' are simply given labels and do not necessarily reflect any geometric notions of dimension. We refer to $n$ as the \emph{filtered dimension} of $X$, or simply as the ``dimension'' when no confusion should arise. 
The set $X^i$ is called the $i$th \emph{skeleton} of $X$, and $X_i=X^i-X^{i-1}$ is the $i$th \emph{stratum}.

\begin{remark}\label{R: filter}
Our definition of a filtered space is more specific than that found in, e.g., \cite{Hug} and other papers by Hughes in that we require $X$ to have a finite number of strata and that the strata be totally ordered by their ``dimensions''. If the skeleta $X^i$ are closed in $X$, then these spaces will also be ``stratified spaces satisfying the frontier condition'' - see \cite{Hug}.  
\end{remark}

A \emph{perversity} $\bar p$ is a function $\bar p: \Z^{\geq 1}\to \Z$ such that $\bar p(k)\leq \bar p(k+1)\leq \bar p(k)+1$. A \emph{traditional perversity} also satisfies $\bar p (1)=\bar p(2)=0$. For  summaries of which topological invariance results on intersection homology hold for which perversities, see \cite{GBF10, Ki}. We shall generally not be concerned with topological invariance in this paper, though we will strongly use that intersection homology with compact supports is a stratum-preserving homotopy invariant; see \cite{GBF3}.

Given $\bar p$ and $X$, one defines $I^{\bar p}C^c_*(X)\subset C^c_*(X)$, the complex of singular chains on $X$, as follows: A simplex $\sigma:\Delta^i\to X$ in $C^c_i(X)$ is \emph{allowable} if $$\sigma^{-1}(X^{n-k}-X^{n-k-1})\subset \{i-k+\bar p(k) \text{ skeleton of } \Delta^i\}.$$ The chain $\xi\in C^c_i(X)$ is allowable if each simplex in $\xi$ and $\bd \xi$ is allowable. $I^{\bar p}C_*^c(X)$ is the complex of allowable chains. $I^{\bar p}C_*^{\infty}(X)$ is defined similarly as the complex of allowable chains in $C_*^{\infty}(X)$, the complex of locally-finite singular chains (also called the complex of chains with \emph{closed support}). Chains in $C_*^{\infty}(X)$ may be composed of an infinite number of simplices (with their coefficients), but for each such chain $\xi$, each point in $X$ must have a neighborhood that intersects only a finite number of simplices (with non-zero coefficients) in $\xi$. See \cite{GBF10} for more details. 

The associated homology theories are denoted $I^{\bar p}H^c_*(X)$ and $I^{\bar p}H_*^{\infty}(X)$. We will sometimes omit the decorations $c$ or $\infty$ if these theories are equivalent, e.g. if $X$ is compact. We will also often omit explicit reference to $\bar p$ below, since all results hold for any fixed perversity. 

Relative intersection homology is defined similarly, though we note that 
\begin{enumerate}
\item the filtration on the subspace will always be that inherited from the larger space by restriction, and
\item  in the closed support case, all chains are required to be locally-finite in the larger space. 
\end{enumerate}

If $(X,A)$ is such a filtered space pair, we use the notation $IC_*^{\infty}(A_X)$ to denote the allowable singular chains supported in $A$ that are locally-finite in $X$. The homology of this complex is $IH_*^{\infty}(A_X)$. Note that in the compact support case, the local-finiteness condition is satisfied automatically so we do not need this notation and may unambiguously refer to $IH_*^c(A)$. The injection  $0\to IC_*^{\infty}(A_X)\to IC_*^{\infty}(X)$ yields a quotient complex $IC_*^{\infty}(X,A)$ and a long exact sequence of intersection homology groups $\to IH_i^{\infty}(A_X)\to IH_i^{\infty}(X)\to IH_i^{\infty}(X,A)\to$. 

If $X$ and $Y$ are two filtered spaces, we call a map $f:X\to Y$
\emph{filtered} if the image of each component of a stratum of $X$
lies in a stratum of $Y$. \emph{N.B. This property is often referred to as ``stratum-preserving'', e.g. in \cite{Q1} and \cite{GBF3}. However, we must reserve  the term ``stratum-preserving'' for other common uses.}
In general, it is not required
that a filtered map take  strata of $X$ to strata of $Y$ of the same (co)dimension. However, if $f$ preserves codimension, or if  $X$ and $Y$ have the same filtered dimension and  $f(X_i)\subset Y_i$, then $f$ will induce a well-defined map on intersection homology (see \cite[Prop. 2.1]{GBF3} for a proof). In this case, we will call $f$ \emph{well-filtered}.
We call a well-filtered map $f$ a \emph{stratum-preserving homotopy equivalence} if there is a well-filtered
map $g:Y\to X$ such that $fg$ and $gf$ are homotopic to the appropriate
identity maps by well-filtered homotopies, supposing that $X\times I$ and $Y\times I$ are given the obvious product filtrations. Stratum-preserving homotopy equivalences induce intersection homology isomorphisms \cite{GBF3}. If stratum-preserving homotopy equivalences between $X$ and $Y$ exist,
we say that $X$ and $Y$ are \emph{stratum-preserving homotopy equivalent}, $X\sim_{sphe} Y$.  

In the sequel, all maps inducing intersection homology homomorphisms will clearly be well-filtered. Hence, we will usually dispense with explicit discussion of this point.  

It is shown in \cite{GBF10} that one can construct a sheaf of intersection chains $\mc{IS^*}$ on any filtered Hausdorff space $X$ such that if $X$ is also paracompact and of finite cohomological dimension then the hypercohomology $\H^*(\mc{IS^*})$ is isomorphic to $IH_{n-*}^{\infty}(X)$, where $n$ is the filtered dimension of $X$. If $X$ is also locally-compact, then $\H^*_c(\mc{IS^*})\cong IH_{n-*}^c(X)$. We will use some properties of these sheaves below, but we refer the reader to \cite{GBF10} for more detailed background.

\subsubsection{A note on coefficients}

Throughout this paper we will leave the coefficient systems tacit so as not to further overburden the notation. However, all results hold for any of the following choices of coefficients, where $R$ is any ring with unit of finite cohomological dimension:

\begin{itemize}
\item Any constant coefficient groups or $R$-modules.

\item Any local system of coefficients of groups or $R$-modules defined on $X-X^{n-1}$ (see \cite{GM2,Bo,GBF10}).

\item If $\bar p$ is a superperversity (i.e. $\bar p(2)>0$; see \cite{CS, GBF10, GBF11}), any stratified system of coefficients $\mc G_0$ as defined in \cite{GBF10} such that $\mc G_0|_{X-X^{n-1}}$ is a local coefficient system of groups or $R$-modules and  $\mc G_0|_{X^{n-1}}=0$. It is shown in \cite{GBF10} that this last coefficient system allows us to recover from singular chains the superperverse sheaf intersection cohomology on pseudomanifolds.

\end{itemize}

\subsection{Stratified homotopies and fibrations}

If $X$ is a filtered space, a map $f:Z\times A\to X$ is \emph{stratum-preserving along $A$} if, for each $z\in Z$, $f(z\times A)$ lies in a single stratum of $X$. If $A=I=[0,1]$, we call $f$ a \emph{stratum-preserving homotopy}. If $f:Z\times I\to X$ is only stratum-preserving when restricted to $Z\times [0,1)$, we say $f$ is \emph{nearly stratum-preserving}.

If $X$ and $Y$ are stratified spaces, a map $p:X\to Y$ is a \emph{stratified fibration} if it admits solutions to stratified lifting problems, i.e. if given a commuting diagram of maps
\begin{equation*}
\begin{CD}
Z&@>f>>& X\\
@V \times 0 VV && @VV p V\\
Z\times I &@>F>>&Y,
\end{CD}
\end{equation*}
such that  $Z$ is any space and $F$ is a stratum-preserving homotopy, there  exists a stratum-preserving homotopy $\td F:Z\times I\to X$ such that $p\td F=F$ and $\td F|_{Z\times 0}=f$. 

 See \cite{Hug, GBF3} for more on stratified fibrations.

\subsection{Manifold homotopically stratified spaces}\label{S: def MHSSs}

Even though the above definition of intersection homology applies to very general spaces, one usually needs to limit oneself to smaller classes of spaces in order to obtain nice properties. Thus at least some of the spaces we will work with will need to satisfy stronger geometric properties, although we don't need the rigidity of, say, the pseudomanifolds on which intersection homology was initially defined  by Goresky and MacPherson \cite{GM1,GM2}. We can often get by with the manifold homotopically stratified spaces introduced by Quinn and refined by Hughes. These spaces were introduced partly with the purpose in mind of being the ``right category'' for intersection homology - see \cite{Q2}.

There is disagreement in the literature as to what to call these spaces. Quinn, himself, calls them both ``manifold homotopically stratified sets'' \cite{Q1} and ``weakly stratified sets''  \cite{Q2}. Hughes \cite{Hu02} prefers the term ``manifold stratified spaces''. We use the term \emph{manifold homotopically stratified space} (MHSS), which seems to capture both that they are stratified by manifolds and that there are additional homotopy conditions on the ``gluing''. 

To define these spaces, we need some preliminary terminology. Except where noted, we take these definitions largely from \cite{Hu02}, with slight modifications to reflect the restrictions mentioned above in Remark \ref{R: filter}:

\subsubsection{Forward tameness and homotopy links}

If $X$ is a filtered space, then $Y$ is \emph{forward tame} in $X$ if there is a neighborhood $U$ of $Y$  and a nearly-stratum preserving deformation retraction $R:U\times I\to X$ retracting $U$ to $Y$ rel $Y$. If the deformation retraction keeps $U$ in $U$, we call $U$ a \emph{nearly stratum-preserving deformation retract neighborhood (NSDRN)}. This last definition was introduced in \cite{GBF5}

The \emph{stratified homotopy link} of $Y$ in $X$ is the space (with compact-open topology) of nearly stratum-preserving paths with their tails in $Y$ and their heads in $X-Y$: $$\hl_s(X,Y)=\{\omega\in X^I\mid \omega(0)\in Y, \omega((0,1])\subset X-Y\}.$$ The \emph{holink evaluation} map takes a path $\omega\in \hl_s(X,Y)$ to $\omega(0)$. For $x\in X_i$, the \emph{local holink}, denoted $\hl_s(X,x)$, is simply the subset of paths $\omega\in \hl_s(X,X_i)$ such that $\omega(0)=x$. Holinks inherit natural stratifications from their defining spaces: $$\hl_s(X,Y)_j=\{\omega\in\hl(X,Y)\mid \omega(1)\in X_j\}.$$

If $X$ is metric and $\delta: Y\to(0,\infty)$ is a continuous function, then $\hl_s^{\delta}(X,Y)$ is the subset of paths $\omega\in\hl_s(X,Y)$ such that $\omega(I)$ is contained inside the open ball $B_{\delta(\omega(0))}(\omega(0))$ of radius $\delta(\omega(0))$ and center $\omega(0)$.

\subsubsection{Manifold Homotopically Stratified Spaces (MHSSs)}
A filtered space $X$ is a \emph{manifold homotopically stratified space (MHSS)} if the following conditions hold:
\begin{itemize}
\item $X$ is locally-compact, separable, and metric.

\item $X$ has finitely many strata, and each $X_i$ is an $i$-manifold without boundary and is locally-closed in $X$. 

\item For each $k>i$, $X_i$ is forward tame in $X_i\cup X_k$. 

\item For each $k>i$, the holink evaluation $\hl_s(X_i\cup X_k,X_i)\to X_i$ is a fibration. 

\item \label{I: CD} For each $x$, there is a stratum-preserving homotopy $\hl(X,x)\times I\to \hl(X,x)$ from the identity into a compact subset of $\hl(X,x)$.  
\end{itemize}

\begin{remark} Condition \eqref{I: CD}, requiring \emph{compactly dominated local holinks}, was not part of the original definition of Quinn \cite{Q1}. It first appears in the work of Hughes leading towards his approximate tubular neighborhood theorem in \cite{Hu02}. While this condition is necessary for  proving the existence of these neighborhoods (when the other dimension conditions of \cite{Hu02} are satisfied), it is not clear that this condition is necessary for any of the other uses of these spaces we make in this paper.
\end{remark}

We say that a MHSS $X$ is $n$-dimensional if its top manifold stratum has dimension $n$. This implies that $X$ is $n$-dimensional in the sense of covering dimension by \cite[Theorem III.2]{HuW}, which states that  a space that is the union of a countable number of closed subsets of dimension $\leq n$ has dimension $\leq n$. This condition holds for $X$ since each stratum is a separable manifold of dimension $\leq n$ (see also \cite[Theorem V.1]{HuW}). It then follows from \cite[Theorem III.1]{HuW} and \cite[Corollary II.16.34, Definition II.16.6, and Proposition II.16.15]{Br} that the cohomological dimension $\dim_{R} X$  of $X$ is $\leq n$ for any ring $R$ with unity (note that since $X$ is metric, it satisfies any paracompactness properties one would like). Similarly, $\dim_R Z\leq n$ for any subspace of $X$. 

A subset of an MHSS is \emph{pure} if it is a closed union of components of strata.

\subsection{Teardrops}

 Given a map $p: X\to Y\times \R$, the \emph{teardrop} $X\cup_p Y$ of $p$ is the space $X\amalg Y$ with the minimal topology such that 
\begin{itemize}
\item $X\into X\cup_p Y$ is an open embedding, and
\item the function $c: X\cup_p Y\to Y\times (-\infty,\infty]$ defined by
\begin{equation*}
c(x)=\begin{cases}
p(x), &x\in X\\
(y,\infty), & y\in Y
\end{cases}
\end{equation*}
\end{itemize}
is continuous.

Given $f:X\to Y$, the teardrop $(X\times \R)\cup_{f\times \text{id}}Y$ is the \emph{open mapping cylinder of $f$ with the teardrop topology}. If $f$ is a proper map between locally compact Hausdorff spaces, then this is the usual mapping cylinder with the quotient topology (see \cite{Hu99a}). An alternative description of the teardrop topology of a mapping cylinder is as the topology  on $X\times (0,1)\amalg Y$  generated by the open subsets of $X\times (0,1)$ and sets of the form $U\cup (p^{-1}(U)\times (0, \epsilon))$, where $U$ is open in $Y$.

If $N$ is a nearly stratum-preserving deformation retract neighborhood (NSDRN) of a pure subset $Y$ of a manifold homotopically stratified space (MHSS), then $N$ is stratum-preserving homotopy equivalent to the mapping cylinder $M$ of the holink evaluation $\hl_s(N,Y)\to Y$, provided $M$ is given the teardrop topology. A proof can be found in \cite[Appendix]{GBF3}.

\section{$IH^c_*$ of stratified fibrations}\label{S: strat fib}

In this section we investigate the compact support singular intersection homology  of a stratified fibration $p:E\to Y$. If $E$ is a particularly nice space, such as a locally-compact homotopically stratified space, then $IH_*^c(E)$ is given by the hypercohomology with compact supports of a sheaf constructed in \cite{GBF10}. Then, in order to find a spectral sequence for $IH_*^c(E)$ starting with $E^2$ groups given in terms of the homology of $Y$ and the intersection homology of the fibers of $p$, we could employ the standard Leray sheaf machinery of, e.g., \cite[Chapter IV]{Br}. However, we will be interested in more general spaces, particularly path spaces. In this case the failure of local compactness on $E$ implies that the arguments of \cite{GBF10} do not hold to show that the hypercohomology with compact supports of the sheaf gives the intersection homology. More specifically, we cannot use Theorem I.6.2 of \cite{Br} to show that $\Gamma_c(\mc{IC}^*)=IC^c_{n-*}(E)$, since $c$ will no longer be a paracompactifying family of supports. Thus we seek an alternative approach by constructing an appropriate sheaf directly  on $Y$. 

We assume throughout that $Y$ is paracompact.\footnote{This is not strictly necessary for everything that follows, but this assumption will hold for all applications we have in mind and will allow us to avoid other, more technical, hypotheses later, such as fine points about the types of rings and supports used to define homological dimensions.}

\subsection{The sheaf and the spectral sequence}

Let $p:E \to Y$ be a stratified fibration to a paracompact space $Y$. Let $n$ be the filtered dimension of $E$.

We define a complex of presheaves $I^{\bar p}_{fc}Z^*$ on $Y$. The definition is in terms of singular intersection chains, but we use $IZ$ instead of $IS$ (the notation of \cite{GBF10}) in order to reserve $IS$ for the intersection chain presheaf defined on $E$.  The decoration $fc$ stands for fiber-wise compact, though the definition is actually a bit more complicated than that.

Let $I^{\bar p}C^{fc}_{n-i}(E)$ be the subgroup of $\xi\in I^{\bar p}C^{\infty}_{n-i}(E)$  satisfying the following support condition: for each point $y\in Y$, there is a neighborhood $N$ of $y$ such that $|\xi|\cap p^{-1}(N)\in c\cap p^{-1}(N)$, the collection of subsets of $E$ that are intersections of $p^{-1}(N)$ with compact sets in $E$. (This definition is motivated by \cite[Chapter IV]{Br}.)  
Similarly, given an open set $U\subset Y$, let $I^{\bar p}C^{fc}_{n-i}((E- \overline{p^{-1}(U)})_E)$ be  the subgroup of $\xi\in I^{\bar p}C^{fc}_{n-i}(E)$ such that $|\xi|\subset E-\overline{p^{-1}(U)}$. 
Let $I^{\bar p}C^{fc}_{n-i}(E,E- \overline{p^{-1}(U)})$ be the quotient   $I^{\bar p}C^{fc}_{n-i}(E)/I^{\bar p}C^{fc}_{n-i}((E- \overline{p^{-1}(U)})_E)$. Then we define the presheaf $I^{\bar p}_{fc}Z^*$ on $Y$ by 

$$U \to I^{\bar p}_{fc}Z^i(U)=I^{\bar p}C^{fc}_{n-i}(E,E- \overline{p^{-1}(U)}).$$

The presheaf restrictions are the obvious ones. Note that for $U\subset V$, the inclusions $I^{\bar p}C^{fc}_{n-i}(U_E) \subset I^{\bar p}C^{fc}_{n-i}(V_E)$ are allowed only because all chains are required to be locally-finite in $E$. Without this restriction, such inclusions would not always be possible, as a locally-finite chain in $U$ might have an accumulation point in $V$.

This presheaf is the most natural one to work with even though our current goal is to study intersection homology with compact supports. This is because the base $Y$ need not be compact, so even if we start with a presheaf defined in terms of compactly supported chains, we would get the same sheaf under sheafification as we do with $I_{fc}Z^*$. To see this,  we also define the complex of presheaves  $I^{\bar p}_cZ^*$ on $Y$ by:
\begin{equation}\label{E: presheaf}
I^{\bar p}_{c}Z^i(U)=I^{\bar p}C^c_{n-i}(E,E- \overline{p^{-1}(U)}).
\end{equation}  
These are defined using the standard compactly supported intersection chains. 

The presheaves $I^{\bar p}_{fc}Z^*$ and $I^{\bar p}_{c}Z^*$ generates sheaves $\mc I^{\bar p}_{fc}\mc Z^*$ and $\mc I^{\bar p}_{c}\mc Z^*$. The inclusion homomorphism $i:I^{\bar p}_{c}Z^*\to I^{\bar p}_{fc}Z^*$ induces a corresponding map of sheaves $i:\mc I^{\bar p}_{c}\mc Z^*\to \mc I^{\bar p}_{fc}\mc Z^*$. This is an isomorphism:

\begin{lemma}\label{L: compact is fc}
The homomorphism $i:\mc I^{\bar p}_{c}\mc Z^*\to \mc I^{\bar p}_{fc}\mc Z^*$ is a sheaf isomorphism.
\end{lemma}
\begin{proof}
The proof is similar to that of \cite[Lemma 3.1]{GBF10}. We must show that $i$ induces an isomorphism at each stalk. 

First, we show injectivity. Let $y\in Y$ and $s\in \mc I^{\bar p}_{c}\mc Z^{n-j}_y$, the stalk at $y$. Suppose that $U$ is a neighborhood of $y$ and  $\xi\in IC^c_j(E, E-\overline{p^{-1}(U)})$ is a finite  chain that represents $s$. If $i|_y(s)=0$, then $\xi=0$ in $IC^{fc}_j(E,E-\overline{p^{-1}(V)})$ for some open $V$ such that $y\in V\subset U$. But this would imply that $|\xi|\subset E-\overline{p^{-1}(V)}$, which implies that $\xi=0$ in $ IC^c_j(E, E-\overline{p^{-1}(V)})$. Hence $s=0$. 

For surjectivity, let $s\in \mc I^{\bar p}_{fc}\mc Z^{n-j}_y$, and suppose $U$ is a neighborhood of $y$ and $\xi\in IC^{fc}_j(E,E-\overline{p^{-1}(U)})$ represents $s$. Due to the support condition, by taking a smaller $U$ if necessary, we may assume that $p^{-1}(U)$ intersects the supports of only a finite number of the simplices of $\xi$. It will suffice to find a finite chain $\zeta\in IC^c_j(E,E-\overline{p^{-1}(U)})$ such that $i(\zeta)=\xi \in IC^{fc}_j(E,E-\overline{p^{-1}(U)})$. If $\xi$ is already a finite chain then $\zeta=\xi$ suffices. Suppose then that $\xi$ contains an infinite number of singular simplices. Let $\Xi$ be the singular chain (not necessarily allowable) composed of singular simplices of $\xi$ (with their coefficients) whose supports intersect $\overline{p^{-1}(U)}$. Let $\xi'$ be the generalized barycentric subdivision of $\xi$ \emph{holding $\Xi$ fixed}. In other words, we perform  a barycentric subdivision of each simplex in $\xi$ except that we do not subdivide the simplices of $\Xi$ nor any common faces between simplices  in $\Xi$ and simplices not in $\Xi$ (see  \cite[\S 16]{MK}). Now take as $\zeta$ the ``regular neighborhood'' of $\Xi$ in $\xi'$. By this we mean  the chain consisting of the simplices in $\Xi$ (with their coefficients) and all other simplices in $\xi'$ that share a vertex  with a simplex in $\Xi$. This $\zeta$ must be finite since $\xi$ is locally finite and $\Xi$ is finite. Furthermore, $\zeta$ is allowable by the arguments in the proof of  \cite[Lemma 2.6]{GBF10}. To see that $i(\zeta)=\xi$ in $IC^{fc}_j(E,E-\overline{p^{-1}(U)})$, we simply note that $\xi-i(\zeta)$ has support in $|\xi-\Xi|$, which lies in $E-\overline{p^{-1}(U)}$. Hence $\xi-i(\zeta)=0$ in $IC^{fc}_j(E,E-\overline{p^{-1}(U)})$.
\end{proof}

The presheaf $IZ_c^*$ with compactly supported chains is the most useful for performing certain local stalk cohomology computations, but the presheaf $IZ_{fc}^*$ with only fiberwise compact supports has some nicer abstract properties, such as the following.

\begin{lemma}\label{L: conj cov}
$I^{\bar p}_{fc}Z^*$ has no non-zero global sections with empty support, and it is conjunctive for coverings.
\end{lemma}
\begin{proof}
The first statement follows similarly to \cite[Lemma 3.2]{GBF10} for intersection chains or \cite[Chapter I, Exercise 12]{Br} for ordinary chains: 
Suppose  $\xi\in I^{\bar p}_{fc}Z^i(Y)=I^{\bar p}C^{fc}_{n-i}(E)$, and suppose that $\xi$ has empty support (meaning the support is empty of the image section under sheafification of the global section of the presheaf represented by $\xi$). So for each point $y\in Y$, the image of $\xi$ in $\lim_{y\in U} IC^{\infty}_{n-*}(E, E-\overline{p^{-1}(U)})=0$. So for each $y\in Y$, there is a neighborhood $U_y$ of $y$ such that the support of $\xi$  lies in $E-\overline{p^{-1}(U_y)}$. In particular, then, $|\xi|\subset \cap_{y\in Y} (Y-\overline{p^{-1}(U_y)})=\emptyset$ and thus $\xi=0$. 

The proof of the second statement is also essentially that of \cite[Lemma 3.3]{GBF10}, which itself is a direct generalization of the usual statement for ordinary singular chains (see Swan \cite[p. 83]{SW}). The present lemma is, in fact,  a special case of  the cited one, which provides the analogous statement for the intersection chain presheaf of $E$ (or any Hausdorff filtered space) by piecing together the bits of chains that are observed in open sets of a covering. The present lemma restricts us to the case in which the open covering of $E$ has the special property that  its sets are of the form $p^{-1}(U)$, $U\subset Y$. The global section we get by piecing together sections over the open sets is a chain in $E$. Thus, it is an element of $I^{\bar p}_{fc}Z^*(Y)$. The only other new concern is the support condition, but this clearly must hold for the constructed global chain, since the support condition is a local one (over $Y$) that is already observed locally in each element of the covering. 
\end{proof}

\begin{corollary}\label{C: agreement}
$\Gamma(Y;\mc I^{\bar p}_{fc}\mc Z^*)=I^{\bar p}C^{fc}_{n-*}(E)$, and if  $Y$ is locally-compact then  $\Gamma_c(Y;\mc I^{\bar p}_{fc}\mc Z^*)=I^{\bar p}C^{c}_{n-*}(E)$.
\end{corollary}
\begin{proof}
The first statement follows from the definition of $I^{\bar p}_{fc} Z^*$ and from \cite[Theorem I.6.2]{Br}, by which the global sections of a presheaf agree with the global sections of the sheaf it induces, provided that it has no nontrivial global sections with empty support and that it is conjunctive for coverings. Note that the collection of closed sets on $Y$ is paracompactifying, since $Y$ is paracompact. Similarly, since $Y$ is locally compact, the family of compact supports is paracompactifying, and so, by the same theorem, $\Gamma_c(Y;\mc I^{\bar p}_{fc}\mc Z^*)$ is equal to $\{s\in I^{\bar p}_{fc} Z^*(Y)=I^{\bar p}C^{fc}_{n-*}(E)\mid |s| \text{ is compact in $Y$}\}$. We need to show that this group is equal to $I^{\bar p}C^{c}_{n-*}(E)$.

Clearly, $I^{\bar p}C^{c}_{n-i}(E)\subset I^{\bar p}C^{fc}_{n-i}(E)$, since if $\xi \in I^{\bar p}C^{c}_{n-i}(E)$, then for each $y\in Y$ and any neighborhood $N$ of $y$, $|\xi|\cap p^{-1}(N)$ is itself of the required form $K\cap p^{-1}(N)$ for a compact set $K$. 

Now suppose that $\xi \in I^{\bar p}C^{fc}_{n-i}(E)$, and let $s$ stand for $\xi$ as an element of  $I^{\bar p}_{fc} Z^*(Y)$ (we make this notational distinction since the geometric support $|\xi|$ of $\xi$ as a chain is a subset of $E$, while the support $|s|$ of $s$ as a presheaf section is a subset of $Y$). Suppose also that $|s|$ is compact. We must show that $|\xi|$ is itself compact, which will complete the proof of the corollary.

Suppose that $y\notin |s|$. Then there is a neighborhood $U$ of $y$ such that $\bar U\cap |s|=\emptyset$ (since $Y$ is paracompact) and such that the image of $\xi$ in $I^{\bar p}C^{fc}_{n-i}(E,E- \overline{p^{-1}(U)})$ is $0$. So $|\xi|\subset E-\overline{p^{-1}(U)}\subset E-p^{-1}(U)$. Thus we see that $|\xi|\subset p^{-1}(|s|)$. For each point $y\in|s|$, there is a neighborhood $N_y$ such that $|\xi|\cap p^{-1}(N_y)=c_y\cap p^{-1}(N_y)$, where $c_y$ is some compact set in $E$. Since $|s|$ is compact, $|s|$ is covered by a finite number of the $N_y$: $|s|\subset \cup_{i=1}^k N_{y_i}$. Thus $|\xi|\subset \cup_{i=1}^k (|\xi|\cap N_{y_i})\subset \cup_{i=1}^k (|\xi|\cap c_{y_i})\subset \cup_{i=1}^k c_i$, which is compact. Thus $|\xi|$ is a closed subset of a compact space and hence compact. 
\end{proof}

\begin{lemma}\label{L: hom fine}
If $E$ is Hausdorff, then $\mc I^{\bar p}_{fc}\mc Z^*$ is homotopically fine.
\end{lemma}
\begin{proof}

This follows essentially as for Proposition 3.5 of \cite{GBF10}, which shows that the intersection chain sheaf on any filtered Hausdorff space $X$ is homotopically fine. In that proof, one begins with the presheaf of compactly supported chains on $X$,  $KS^*(U)=IC^c_{n-*}(X,X-\bar U)$, and shows that each singular chain is chain homotopic to a sum of chains $\mc U$ small for a locally-finite covering $\mc U$ of $X$. One then verifies that this chain homotopy persists at the sheaf level.

In the case under consideration, although our base space is $Y$, our chains live in $E$, and  the presheaf $I_{c}^{\bar p}Z^*$ on $Y$ is essentially the presheaf $KS^*$ on $E$, with the allowable open set inputs restricted to be of the form $p^{-1}(U)$ for $U$ open in $Y$. Thus the arguments of Proposition 3.5 of  \cite{GBF10} mostly go through in the same fashion as a special case. The sheafification is slightly different, but the same arguments carry through with minor modifications.
\end{proof}

Thus we get the following theorem.

\begin{theorem}\label{T: strat fib SS}
If $E\to Y$ is a stratified fibration with $E$ Hausdorff  and $Y$ paracompact of finite cohomological dimension, then $I^{\bar p}H^{fc}_{n-*}(E)$ is the abutment of a spectral sequence with $E^2$ terms $H^p(Y;\mc{H}^q(\mc I^{\bar p}_{fc}\mc Z^*))$. If $Y$ is also locally compact, then $I^{\bar p}H^c_{n-*}(E)$ is  the abutment of a  spectral sequence with $E^2$ terms $H^p_c(Y;\mc{H}^q(\mc I^{\bar p}_{fc}\mc Z^*))$.
\end{theorem}
\begin{proof}
Since $\mc I^{\bar p}_{fc}\mc Z^*$ is homotopically fine, it follows that $H^*(H^p(Y; \mc I^{\bar p}_{fc}\mc Z^*))=H_c^*(H^p(Y; \mc I^{\bar p}_{fc}\mc Z^*))=0$ for all $p>0$ by \cite[p. 172]{Br}, provided the respective families of supports are paracompactifying. But this is ensured by the hypotheses in the respective cases. So if $Y$ has finite dimension, then the hypercohomology spectral sequences 
with $E_2^{p,q}=H^p(Y;\mc{H}^q(\mc I^{\bar p}_{fc}\mc Z^*))$ and $E_2^{p,q}=H_c^p(Y;\mc{H}^q(\mc I^{\bar p}_{fc}\mc Z^*))$ respectively abut to $H^{p+q}(\Gamma(Y;\mc I^{\bar p}_{fc}\mc Z^*))$ and $H^{p+q}(\Gamma_c(Y;\mc I^{\bar p}_{fc}\mc Z^*))$ by \cite[IV.2.1]{Br}. But by Corollary \ref{C: agreement}, these groups are respectively $I^{\bar p}H^{fc}_{n-p-q}(E)$ and $I^{\bar p}H^{c}_{n-p-q}(E)$.
\end{proof}

\subsection{The stalk cohomology}\label{S: stalks}

In order to employ Theorem \ref{T: strat fib SS} usefully, it is necessary to understand the stalks of $\mc I^{\bar p}_{fc}\mc Z^*$ and their cohomology.

\begin{lemma}\label{L: C stalk}
For $y\in Y$, the stalk $\mc I^{\bar p}_{fc}\mc Z^*_y$ is isomorphic to $IC^c_{n-*}(E,E-p^{-1}(y))$.
\end{lemma}
\begin{proof}
By Lemma \ref{L: compact is fc}, $\mc I^{\bar p}_{fc}\mc Z^*_y$  is $\lim_{y\in U} I^{\bar p}C^c_{n-*}(E, E-\overline{p^{-1}(U)})$. Since direct limits are exact as functors, this is the quotient of $I^{\bar p}C^c_{n-*}(E)$ by $\lim_{y\in U} I^{\bar p}C^c_{n-*}(E-\overline{p^{-1}(U)})$, which, as for ordinary singular chains, is $I^{\bar p}C^c_{n-*}(\cup_U E-\overline{p^{-1}(U)})=I^{\bar p}C^c_{n-*}(E-p^{-1}(y))$.
\end{proof}

\begin{corollary}\label{C: H stalk}
The stalk  at $y$ of the derived sheaf $\mc H^*(\mc I^{\bar p}_{fc}\mc Z^*)$ is isomorphic to $IH^c_{n-*}(E,E-p^{-1}(y))$.
\end{corollary}

We next show that if the base $Y$ of the stratified fibration $p: E\to Y$ is a manifold homotopically stratified space, then the sheaf $\mc I^{\bar p}_{fc}\mc Z^*$ on $Y$ is cohomologically locally constant on each stratum of $Y$. Note that we must require some such nice conditions on the base in order to obtain this uniformity. In ordinary fibration theory, one might only need local contractibility of the base in order to get homotopy local triviality of the bundle. In the stratified case, we don't have stratified contractibility (unless the base is actually unstratified), so we must impose some other stratified uniformity property. It would be interesting to know the weakest conditions under which our following results  hold, but assuming $Y$ to be an MHSS will be suitable for our purposes and seems reasonable for studying intersection homology. Essentially, we only use that the strata are manifolds and that we have a way of extending homotopies given on strata.

\begin{lemma}\label{L: loc con}
Let $Y$ be an MHSS. Let $(B,A)\subset Y_k$ be a closed subset pair homeomorphic to $(D^k, \frac{1}{2}D^k)$, where $D^k$ is the $k$-dimensional ball. Let  $p:E\to Y$ be a stratified fibration. Then for each $y\in \text{int}(A)$, restriction $IH^c_*(E, E-p^{-1}(A))\to IH^c_*(E,E-p^{-1}(y))$ is an isomorphism.
\end{lemma}
\begin{proof}
We will employ the long exact sequence of the triple and show that $IH^c_*(E-p^{-1}(y),E-p^{-1}(A))=0$.  We need only show that the support of any relative cycle $\xi$ in $IC^c_*(E-p^{-1}(y),E-p^{-1}(A))$ can be stratum-preserving homotoped into $p^{-1}(Y)-p^{-1}(A)$. This homotopy can then be used to establish a homology as in \cite{GBF3}.

So let $[\xi]$ be a class in $IH^c_*(E-p^{-1}(y),E-p^{-1}(A))$. Certainly there is a deformation of $Y_k$ that takes $p(|\xi|)\cap A$ to $Y-A$ and such that the trace of the homotopy does not contain $y$ - just push $p(|\xi|)\cap Y_k$ ``radially'' away from $y$ in $B$ by an isotopy of $B$ rel boundary. Considering this deformation as a stratum-preserving deformation of $Y^k$ (just hold everything outside  $B$ fixed), we can apply \cite[Corollary 6.4]{Hug} to extend the deformation to all of $Y$. We  then lift this homotopy to a stratum-preserving homotopy of $E$ that takes $|\xi|$ to $E-p^{-1}(A)$. 
\end{proof}

\begin{corollary}\label{C: clc}
If $Y$ is an MHSS, then $\mc I^{\bar p}_{fc}\mc Z^*$ is cohomologically locally constant on each stratum of $Y$. 
\end{corollary}
\begin{proof}
Since each stratum  $Y_k$ is a manifold, each point $y_0\in Y_k$ has a neighborhood pair $(B,A)$ in $Y_k$ homeomorphic to $(D^k, \frac{1}{2}D^k)$. By Lemma \ref{L: loc con} and Corollary \ref{C: H stalk}, for each $y\in \text{int}(A)$, 
$\mc H(\mc I^{\bar p}_{fc}\mc Z^*)_y\cong IH^c_{n-*}(E,E-p^{-1}(y))\cong IH^c_{n-*}(E, E-p^{-1}(A))$. It follows that $\mc H(\mc I^{\bar p}_{fc}\mc Z^*)$ is cohomologically locally-constant over $Y_k$.
\end{proof}

\subsection{Relative $IH^c_*$ of stratified fibrations}\label{S: relative}

We also want to study $IH_*(E,A)$, where $(E,A)$ is a stratified fibration pair.

\begin{definition}\label{D: strat fib pair}
$p: (E,A)\to Y$ is a \emph{stratified fibration} pair if it admits solutions to relative stratified lifting problems. In other words, given a relative lifting problem
\begin{diagram}
(Z,X) &\rTo^f & (E,A)\\
\dTo_{\times 0} &&\dTo^p\\
(Z,X)\times I& \rTo^F &Y 
\end{diagram}
in which $F$ is stratum-preserving along $I$ and $f$ is a map of pairs, there exists a map $\td F: (Z,X)\times I \to (E,A)$ that is stratum-preserving along $I$, satisfies $p\td F=F$ and $\td F(z,0)=f(z)$, and takes $X\times I$ into $A$.  Note that this makes $p:A\to Y$ a stratified fibration in its own right.
\end{definition}

To study the intersection homology of such pairs, one employs quotient presheaves. Let $I^{\bar p}_cZ^*_E$ and $I^{\bar p}_cZ^*_A$ be the presheaves as defined by equation \eqref{E: presheaf} with respect to the stratified fibrations $p: E\to Y$ and $p|_A: A\to Y$. One easily checks that there is a well-defined injection $I^{\bar p}_cZ^*_{A}\into I^{\bar p}_cZ^*_{E}$, which induces the quotient presheaf $I^{\bar p}_cZ^*_{E,A}$: $$U\to I^{\bar p}C^c_{n-*}(E,E-\overline{p^{-1}(U)})/ I^{\bar p}C^c_{n-*}(A,A-A\cap \overline{p^{-1}(U)}).$$ 

Our most useful homological tool will be the following lemma, the proof of which employs the techniques of \cite{GBF10}.

\begin{lemma}\label{L: rel presheaf}
For $U$ an open subset of $Y$, $H^*(I^{\bar p}_cZ^*_{E,A}(U))\cong IH^c_{n-*}(E,A \cup (E-\overline{p^{-1}(U)}))$. In particular, $H^*(I^{\bar p}_cZ^*_{E,A}(Y))\cong IH^c_{n-*}(E,A)$. 
\end{lemma}
\begin{proof}
By definition, $I^{\bar p}_cZ^*_{E,A}(U)=I^{\bar p}C^c_{n-*}(E,E- \overline{p^{-1}(U)})/ I^{\bar p}C^c_{n-*}(A,A- A\cap\overline{p^{-1}(U)})$, which is isomorphic to $$\frac{I^{\bar p}C^c_{n-*}(E)}{I^{\bar p}C^c_{n-*}(A) + I^{\bar p}C^c_{n-*}(E-\overline{p^{-1}(U)})}.$$
Using \cite[Corollary 2.10]{GBF10}, the ``denominator'' is chain homotopy equivalent to $IC^c_{n-*}(A \cup (E-\overline{p^{-1}(U)})).$ The lemma now follows from long exact sequences and the five lemma.
\end{proof}

The short exact sequence of presheaves 
\begin{equation*}
\begin{CD}
0 &@>>>& I^{\bar p}_cZ^*_{A} &@>>>&  I^{\bar p}_cZ^*_{E} &@>>>&  I^{\bar p}_cZ^*_{E,A} &@>>>& 0
\end{CD}
\end{equation*}
passes to a short exact sequence of sheaves 
\begin{equation}\label{E: sheaf SES}
\begin{CD}
0 &@>>>& \mc I^{\bar p}_{fc}\mc Z^*_{A} &@>>>&  \mc I^{\bar p}_{fc}\mc Z^*_E &@>>>&  \mc I^{\bar p}_{fc}\mc Z^*_{E,A} &@>>>& 0,
\end{CD}
\end{equation}
in which the first two sheaves are those we have met before (with respect to the fibrations $A\to Y$ and $E\to Y$) and the third is their quotient sheaf.

Note that we won't be able to employ the same arguments we used with $\mc I_{fc}\mc Z^*_E$ to conclude that $\Gamma(Y, \mc I^{\bar p}_{fc}\mc Z^*_{E,A})\cong I^{\bar p}_cZ^*_{E,A}(Y)$, since the relative presheaf is \emph{not} conjunctive for coverings\footnotemark. However, we will obtain what we want by other means in a moment. First, we will see that the hypercohomology is the same as the section cohomology via the following lemma.

\footnotetext{Consider $Y=U\cup V$, and consider an intersection chain consisting of a single simplex $\sigma$ whose support intersects $p^{-1}(U)-p^{-1}(U)\cap A$ and $p^{-1}(V)-p^{-1}(V)\cap A$, but not $p^{-1}(U\cap V)-p^{-1}(U\cap V)\cap A$ (so over $U\cap V$, $\sigma$ lies in $A$). Then the representative of $\sigma$ in $I_cZ^*(U)$ and the section $0$ in $I_cZ^*(V)$ restrict to the same section $0$ of $I_cZ^*(U\cap V)$. But  these sections are not the restrictions of a single element of $I_cZ^*(Y)$.}

\begin{lemma}
$\mc I^{\bar p}_{fc}\mc Z^*_{E,A}$ is homotopically fine.
\end{lemma}
\begin{proof}
By Lemma \ref{L: hom fine}, we already know that $\mc I^{\bar p}_{fc}\mc Z^*_{E}$ and $\mc I^{\bar p}_{fc}\mc Z^*_{A}$ are homotopically fine. We noted in the proof of that lemma that this is essentially a special case of Proposition 3.5 of \cite{GBF10}. That proposition is proven by first showing that, for any Hausdorff filtered space $Z$ with locally-finite open covering $\mc U=\{U_k\}$, the identity map of $IC^c_{n-*}(Z)$ is chain homotopic, by a chain homotopy $D$, to $\sum_k g_k$, where each $g_k$ is an endomorphism of $IC^c_{n-*}(Z)$ with support in $U_k$. The maps $g_k$ and $D$ descend to the sheaf level, providing a chain homotopy to the identity and demonstrating the desired property of the sheaf complex. In the present case, we once again specialize to covers of the form $\{p^{-1}(U_k)\}$, $U_k\subset Y$, as in Lemma \ref{L: hom fine}, and we need only observe that since $|g_k(\xi)|\subset |\xi|$ and $|D(\xi)|\subset |\xi|$, for $g_k$ and $D$ as defined in \cite[Proposition 3.5]{GBF10}, $g_k$ and $D$ each induce endomorphisms of the presheaf $I^{\bar p}_cZ^*_{E,A}$. Once again, these endomorphisms descend to the sheaf level to provide a chain homotopy from the identity on $\mc I^{\bar p}_{fc}\mc Z^*_{E,A}$ to a sum of sheaf endomorphisms with $\mc U$ small support.
\end{proof}

\begin{theorem}\label{T: rel identification}
The long exact hypercohomology sequence
\begin{equation}\label{E: hyper LES}
\begin{CD}
@>>> & \H^i(\mc I^{\bar p}_{fc}\mc Z^*_{A}) &@>>>&  \H^i(\mc I^{\bar p}_{fc}\mc Z^*_E) &@>>>&  \H^i(\mc I^{\bar p}_{fc}\mc Z^*_{E,A}) &@>>>& \H^{i+1}(\mc I^{\bar p}_{fc}\mc Z^*_{A}) )&@>>>
\end{CD}
\end{equation}
is isomorphic to the long exact sequence 
\begin{equation}\label{E: IH LES}
\begin{CD}
@>>> & IH^{fc}_{n-i}(A_E) &@>>>&  IH^{fc}_{n-i}(E) &@>>>& IH^{fc}_{n-i}(E,A) &@>>>& IH^{fc}_{n-i-1}(A) &@>>>.
\end{CD}
\end{equation}
If $Y$ is locally compact and in the first sequence we take hypercohomology with compact supports, then the isomorphism is to the version of the second sequence with compact supports.
\end{theorem}

\begin{proof}
By definition,   \eqref{E: IH LES} is simply the long exact sequence corresponding to the short exact sequence
 \begin{equation*}
\begin{CD}
0 &@>>>& I^{\bar p}C^{fc}_{n-*}(A) &@>>>&  I^{\bar p}C^{fc}_{n-*}(E) &@>>>&  I^{\bar p}C^{fc}_{n-*}(E,A)=I^{\bar p}C^{fc}_{n-*}(E)/I^{\bar p}C^{fc}_{n-*}(A) &@>>>& 0.
\end{CD}
\end{equation*}
Note that the first (non-trivial) map is well-defined since if $\xi$ is a chain in $I^{\bar p}C^{fc}_{n-*}(A)$, then it is certainly allowable in $I^{\bar p}C^{fc}_{n-*}(E)$; the only issue is that $\xi$ be locally-finite in $E$. But let $x$ be a point in the $E$, and let $y=p(x)$. Then by definition of 
$I^{\bar p}C^{fc}_{n-*}(A)$, there is a neighborhood $N$ of $y$ such that $|\xi|\cap p^{-1}(N)\cap A=c\cap p^{-1}(N)$ for some compact set $c\subset A$. But as $\xi$ is also assumed to be locally finite in $A$, this implies that $p^{-1}(N)\cap A$ intersects only a finite number of simplices of $\xi$. Then since $|\xi|\subset A$, $p^{-1}(N)\subset E$ is a neighborhood of $x$ that intersections only a finite number of simplices of $\xi$. Thus $\xi$ is also locally finite in $E$.

Consider now the following map of short exact sequences:
\begin{diagram}
0 &\rTo& I^{\bar p}C^{fc}_{n-*}(A) &\rTo&  I^{\bar p}C^{fc}_{n-*}(E) &\rTo&  I^{\bar p}C^{fc}_{n-*}(E,A)&\rTo& 0\\
&&\dTo&&\dTo&&\dTo\\
0 &\rTo& \Gamma(Y, \mc I^{\bar p}_{fc}\mc Z^*_{A} ) &\rTo& \Gamma(Y, \mc I^{\bar p}_{fc}\mc Z^*_{E} )  &\rTo& Im(\Gamma(Y, \mc I^{\bar p}_{fc}\mc Z^*_{E} ) \to \Gamma(Y, \mc I^{\bar p}_{fc}\mc Z^*_{E,A} ) ) &\rTo& 0.
\end{diagram}
The left two vertical maps are induced by mapping global presheaf sections to global sheaf sections. By Corollary \ref{C: agreement} they are isomorphisms. The righthand vertical map is induced by the others, the bottom sequence being exact because  $\Gamma$ is a left exact functor (as applied here to the sequence \eqref{E: sheaf SES}). The righthand map is also an isomorphism by the five lemma. This isomorphism of short exact sequences induces an isomorphism of the corresponding long exact sequences in homology. So we need to show that the long exact homology sequence of the bottom  short exact sequence is the hypercohomology exact sequence \eqref{E: hyper LES}. 

By \cite[p. 276, \#14, and p. 72]{Br}, the long exact sequence associated to the bottom short exact sequences is isomorphic to the exact sequence  
\begin{diagram}
\rTo& H^i(\Gamma(Y, \mc I^{\bar p}_{fc}\mc Z^*_{A} )) &\rTo& H^i(\Gamma(Y, \mc I^{\bar p}_{fc}\mc Z^*_{E} ) ) &\rTo& H^i(\Gamma(Y, \mc I^{\bar p}_{fc}\mc Z^*_{E,A} ) ) &\rTo& H^{i+1}(\Gamma(Y, \mc I^{\bar p}_{fc}\mc Z^*_{A} )) &\rTo, 
\end{diagram}
 since these sheaves are all homotopically fine. For the same reason, the groups in this sequence are the  hypercohomology groups of the respective sheaves by \cite[p. 172]{Br} and \cite[IV.2.1]{Br} (see also the discussion following Proposition 3.5 of \cite{GBF10}).

The corresponding statements with compact supports hold by the same arguments using Corollary \ref{C: agreement} and the same references in \cite{Br}.
\end{proof}

So we have the following relative analogue of Theorem \ref{T: strat fib SS}.

\begin{theorem}\label{T: rel strat fib SS}
If $(E,A)\to Y$ is a stratified fibration pair with $E$ and $A$ Hausdorff and $Y$ paracompact of finite cohomological dimension, then $I^{\bar p}H^{fc}_{n-*}(E,A)$ is the abutment of a spectral sequence with $E^2$ terms $H^p(Y;\mc{H}^q(\mc I^{\bar p}_{fc}\mc Z_{E,A}^*))$. If $Y$ is also locally compact, then 
$I^{\bar p}H^c_{n-*}(E,A)$ is the abutment of  a spectral sequence with $E^2$ terms  $H^p_c(Y;\mc{H}^q(\mc I^{\bar p}_{fc}\mc Z_{E,A}^*))$.
\end{theorem}
\begin{proof}
Using Theorem \ref{T: rel identification}, this theorem now follows as for Theorem \ref{T: strat fib SS} from the hypercohomology exact sequence; see \cite[IV.2.1]{Br}.
\end{proof}

Next, of course, we want to compute the stalk cohomology of $\mc I^{\bar p}_{fc}\mc Z_{E,A}^*$. Unfortunately, in the relative case the analogue of Lemma \ref{L: C stalk} will be a bit messy. Using that lemma and exactness of direct limits, we see that $$(\mc I^{\bar p}_{fc}\mc Z^*_{E,A})_y\cong \frac{IC^c_{n-*}(E,E-p^{-1}(y))}{IC^c_{n-*}(A, A-A\cap p^{-1}(y))}\cong \frac{IC^c_{n-*}(E)}{IC^c_{n-*}(A)+IC^c_{n-*}(E- p^{-1}(y))}.$$ However, the stalk cohomology  will be more reasonable:

\begin{proposition}\label{P: rel stalk}
The stalk  at $y$ of the derived sheaf $\mc H^*(\mc I^{\bar p}_{fc}\mc Z^*_{E,A})$ is isomorphic to $IH^c_{n-*}(E,A\cup (E-p^{-1}(y)))$.
\end{proposition}
\begin{proof}
This follows by using direct limit arguments as in the proof of Lemma \ref{L: C stalk} together with the homology computations of Lemma \ref{L: rel presheaf}.
\end{proof}

\begin{proposition}\label{P: rel clc}
If $Y$ is an MHSS, then $\mc I^{\bar p}_{fc}\mc Z_{E,A}^*$ is cohomologically locally constant on each stratum of $Y$. 
\end{proposition}
\begin{proof}
This follows  by an argument similar to the proof of Corollary \ref{C: clc} using Lemma \ref{L: rel presheaf}. 
\end{proof}

\begin{corollary}
If $Y$ is an MHSS, then the stalk cohomology sequence at $y$ of the short exact sequence
\begin{equation*}
\begin{CD}
0 &@>>>& \mc I^{\bar p}_{fc}\mc Z^*_{A} &@>>>&  \mc I^{\bar p}_{fc}\mc Z^*_E &@>>>&  \mc I^{\bar p}_{fc}\mc Z^*_{E,A} &@>>>& 0
\end{CD}
\end{equation*}
of stratified cohomologically locally constant sheaves is isomorphic to the long exact sequence
\begin{equation*}
\begin{CD}
 @>>> & IH^c_{n-*}(A,A-A\cap p^{-1}(y))   &@>>> & IH^c_{n-*}(E,E-p^{-1}(y))   &@>>> &  IH^c_{n-*}(E,A\cup (E-p^{-1}(y))) &@>>> 
\end{CD}
\end{equation*}
\end{corollary}

\section{$IH^{\infty}_*$ of stratified fibrations}\label{S: strat fib closed}

In this section we study the closed support intersection homology of stratified fibrations (which is more akin to the Deligne sheaf intersection cohomology of \cite{GM2} and \cite{Bo}). In this case,  we do not a priori have to be so concerned about the possible failure of local compactness of $E$, which led us to construct $IZ^*$ on $Y$ in Section \ref{S: strat fib}. This is because the family of closed supports will be paracompactifying so long as $E$ is paracompact, which allows us to employ more sheaf theoretic machinery  directly on $E$. The sheaves we have so far had to construct directly on $Y$ can instead be taken simply as the direct images of the sheaves on $E$, and so we can employ more standard technology. The main problem, however, is that closed support chains and homotopies (even stratified homotopies) don't get along so well. So uniformity results along the lines of Corollary \ref{C: clc} and Proposition \ref{P: rel clc} will require some extra conditions on $E$ or the stratified fibration $p$.

We let $\mc {IS}^*$ denote the intersection chain sheaf on $E$ as defined in \cite{GBF10}, and we assume $E$ is paracompact of finite cohomological dimension. This sheaf is homotopically fine, and its hypercohomology computes $IH^{\infty}_{n-*}(E)$, the intersection homology of $E$ with locally-finite chains. $\mc{IS}^*$ is the sheafification of both the presheaves $U\to IC^{\infty}_{n-*}(X, X-\bar U)$ and $U\to IC^c_{n-*}(E, E-\bar U)$.

Since $\mc{IS}^*$ is homotopically fine, $IH^{\infty}_{n-*}(E)\cong \H^*(E;\mc{IS}^*)$, the hypercohomology of $\mc{IS^*}$ - see \cite{GBF10} or, in \cite{Br}, the paragraph preceding Theorem IV.2.1 and Exercise 32, page 17. Since $E$ has finite cohomological dimension, injective resolutions exist, and these groups are also isomorphic to  $\H^*(Y;Rp_*\mc{IS}^*)$. This gives us a spectral sequence immediately from the hypercohomology spectral sequence:

\begin{theorem}\label{T: lf IH SS}
If $E\to Y$ is a stratified fibration with $E$  and $Y$ paracompact of finite cohomological dimension, then $I^{\bar p}H^{\infty}_*(E)$ is the abutment of a spectral sequence with $E^2$ terms $H^p(Y;\mc{H}^q( Rp_*\mc{IS}^* ))$.
\end{theorem}

So, we must examine the cohomology stalks of  $Rp_*\mc{IS}^* $.

We have 
\begin{align*}
H^*(Rp_*\mc{IS}^*)_y &=\lim_{y\in U}H^*(\Gamma(U;Rp_*\mc{IS}^*))\\
&=\lim_{y\in U}H^*(\Gamma(p^{-1}(U);\mc J^*)), \text{  where $\mc J^*$ is an injective resolution of $\mc {IS}^*$}\\
&= \lim_{y\in U}\H^*(p^{-1}(U);\mc{IS}^*).
\end{align*}
Since, for an open set $V$, $\mc{IS}^*|V$ is the intersection chain sheaf on $V$ by Proposition 3.7 of \cite{GBF10}, this last group is $\lim_{y\in U} IH^{\infty}_{n-*}(p^{-1}(U))$.

In general, it will be hard to say what this is - stratified fibrations will not necessarily have nice structures over small neighborhoods in the base unless, for example, the base is unstratified (more on this below). Even then, the nice local structures  will generally only be up to stratified homotopy, and homology with  closed supports does not behave well under such homotopies. Below, we will show how we can compute these groups if we are willing to place stronger assumptions on the spaces involved.

\vskip1cm

One case in which we can compute locally-finite intersection homology effectively is that for stratified fibration pairs with cocompact fibers: As in Section \ref{S: relative}, let $p:(E,A)\to Y$ be a stratified fibration pair, and suppose that for each compact $K\subset Y$, $p^{-1}(K)\cap A$ is cocompact in $p^{-1}(K)$, i.e. $p^{-1}(K)-p^{-1}(K)\cap A$ is compact. We  assume that $Y$ is paracompact, locally compact,  and of finite cohomological dimension. In this case, we claim that $IH^{\infty}_*(E,A)\cong IH_*^{fc}(E,A)$, which implies that to obtain $IH^{\infty}_*(E,A)$ we can use the spectral sequence of Theorem \ref{T: rel strat fib SS} and the stalk calculations of Proposition \ref{P: rel stalk}. 

Note that here, as always, relative locally finite intersection homology is defined by
$IH^{\infty}_*(E,A)\cong H_*(IC_*^{\infty}(E)/IC_*^{\infty}(A_E))$, where $IC_*^{\infty}(A_E)$ is the subgroup of $\xi\in IC_*^{\infty}(E)$ such that $|\xi|\subset A$. In particular all chains in $IC_*^{\infty}(A_E)$ are locally finite in $E$.

\begin{proposition}\label{P: cocompact}
Let $p:(E,A)\to Y$ be a stratified fibration pair with $E$, $A$, and $Y$ paracompact of finite cohomological dimension and so that $Y$ is also locally-compact. Suppose also that for each compact $K\subset Y$, $p^{-1}(K)\cap A$ is cocompact in $p^{-1}(K)$. Then $IH^{\infty}_*(E,A)\cong IH_*^{fc}(E,A)$, 
\end{proposition}
\begin{proof}
Recall that $IH_*^{fc}(E,A)$ is defined as $H_*(IC_*^{fc}(E)/IC_*^{fc}(A))$, where the chains in each group of the quotient satisfy the defining support condition. We have compatible homomorphisms (in fact, inclusions) $IC_*^{fc}(E)\to IC_*^{\infty}(E)$ and $IC_*^{fc}(A)\to IC_*^{\infty}(A_E)$ (see the proof of Theorem \ref{T: rel identification})
inducing  a homomorphism $h: IH_*^{fc}(E,A)\to IH^{\infty}_*(E,A)$. We need to show that $h$ is an isomorphism.
The technique is essentially that used in \cite{GBF10} to show that if $Z$ is a filtered space and $X$ is cocompact in $Z$ then $IH^{\infty}_*(Z,X)\cong IH^{c}_*(Z,X)$ - one subdivides finely and then throws out irrelevant cells away from $Z-X$, but we must take some extra care in case $E-A$ is not compact.

We first show that $h$ is surjective. Let $\xi$ be a cycle in $ IC_i^{\infty}(E,A)$.

Since $y$ is locally compact, let us assign to each point $y\in Y$ an open neighborhood $U_y$ such that $\bar U_y$ is compact. Since $Y$ is paracompact, this covering has a locally finite refinement, say $\mc V$. Now, as in the methods of \cite{GBF10}, we can find a subdivision of $\xi$ such that each simplex $\sigma$ of the subdivision satisfies $|\sigma|\subset p^{-1}(V)$ for some $V\in\mc V$. This can be done inductively over skeleta of $\xi$: Clearly it is satisfied for the vertices of simplices in $\xi$. Then suppose that we have subdivided the $j$ skeleta of the simplices in $\xi$ so that each simplex in the subdivision of the $j$ skeleta is mapped into some $p^{-1}(V)$. Then for each $j+1$ simplex $\tau$ of a simplex of $\xi$ we perform a generalized subdivision of $\tau$ holding $\bd \tau$ fixed and such that each $j+1$ subsimplex of the subdivision of $\tau$ is contained in some $p^{-1}(V)$. By induction, we obtain the desired subdivision $\xi'$ of $\xi$. More details can be found in the proof of \cite[Proposition 2.9]{GBF10}. By \cite[Proposition 2.7]{GBF10}, $\xi'$ and $\xi$ are intersection homologous, i.e. $[\xi]=[\xi']\in IH_i^{\infty}(E,A)$. 

Next, let $\Xi$ be the chain consisting of all simplices in $\xi'$ that intersect $E-A$ (with their coefficients from $\xi'$); $\Xi$ need not be allowable. Let $\xi''$ be a barycentric subdivision of $\xi'$ rel $\Xi$, and let $\eta$ be the chain consisting of all simplices of $\xi''$ that share a vertex with a simplex of $\Xi$ (with their coefficients). It follows as in the proof of \cite[Lemma 2.12]{GBF10} that $\eta$ and $\xi''-\eta$ are allowable and so, since $|\xi''-\eta|\subset A$,  $[\eta]=[\xi'']=[\xi]$ (again using \cite[Proposition 2.7]{GBF10} for the second equality). We claim that $\eta\in  IC_*^{fc}(E)$ and hence the class it represents in $IH_*^{fc}(E,A)$ maps to $[\xi]$ under $h$.

So let $y$ be a point in $Y$. Then $y\in W$ for some $W\in \mc V$. Let $\mc V_1=\{V\in \mc V\mid \bar V\cap \bar W\neq \emptyset\}$, and let $W_1=\cup_{V\in \mc V_1} V$. Since each $V\subset U_z$ for some $z\in Y$, $\bar W$ is compact, and thus $\mc V_1$ is a finite collection of sets and $\bar W_1$ is compact. Let $\mc V_2=\{V\in \mc V\mid \bar V\cap \bar W_1\neq \emptyset\}$, and let $W_2=\cup_{V\in \mc V_2} V$. Again, $\mc V_2$ is finite and and $\bar W_2$ is compact. We note that for each $V\in \mc V_2-\mc V_1$, $\bar V\cap \bar W=\emptyset$. 

Now, consider $p^{-1}(W)$. We ask which simplices of $\eta$ can intersect $p^{-1}(W)$. Recall that each simplex of $\eta$ must have its support in some $p^{-1}(V)$. So first there are the simplices of $\Xi$ which intersect $p^{-1}(W)$. By construction of $\xi'$ and $\Xi$, these consist of the simplices that actually lie in  $p^{-1}(W)$ and intersect $E-A$ and those that lie in $p^{-1}(V)$ and intersect $E-A$ for some $V$ such that $V\cap W$ is non-empty. Such $V$ are in $\mc V_1$. Next, there are the simplices $\sigma$ of $\eta$ that intersect $p^{-1}(W)$ and are not in $\Xi$ but share a vertex with a chain $\tau$ in $\Xi$. By the same argument, $|\sigma|$ must be in $p^{-1}(V)$ for some $V\in\mc V_1$, and so the simplex $\tau$ in $\Xi$ that shares a vertex with $\sigma$ must have its support contained in $p^{-1}(V)$ for some $V\in \mc V_2$. Now, since $\bar W_2$ is compact, $p^{-1}(\bar W_2)-p^{-1}(\bar W_2)\cap A$ is compact and the collection $\mf S$ of simplices of $\eta$ whose supports intersect this set is finite. Then, by local finiteness of $\eta$, the collection of simplices $\mf S_1$ in $\eta$ that share vertices with simplices in $\mf S$ is finite. But we have seen that this last collection of simplices contains all simplices that intersection $p^{-1}(W)$.  Thus $\cup_{\sigma\in \mf S_1}|\sigma|$ is a compact set whose intersection with $p^{-1}(W)$ is  $|\eta|\cap p^{-1}(W)$, which shows that $\eta$ has the support structure desired.

The proof that $h$ is injective is similar: given a chain $\xi$ in $IC^{fc}_*(E,A)$ that relatively bounds a chain $\zeta$ in $IC^{\infty}_*(E,A)$, we can use the above procedure to subdivide and truncate $\zeta$ to obtain a chain in $IC^{fc}_*(E,A)$  whose relative boundary is relatively homologous to $\xi$.
\end{proof}

Other special cases in which we can compute $IH^{\infty}_*$ are considered in the following sections.

\section{Manifold bases}\label{S: man base}

In this section, we consider the special case in which the base of a stratified fibration is an unstratified topological manifold. In this case, some of our earlier computations can be stated in a cleaner form.

In particular, we have the following version of Corollary \ref{C: H stalk}.

\begin{proposition}
Let $p: (E,A)\to Y$ be a stratified fibration pair, where $Y$ is an unfiltered $m$-manifold.  Then the stalk at $y$ of the locally constant derived sheaf $\mc H^*(\mc I^{\bar p}_{fc}\mc Z_E^*)$ is isomorphic to $IH^c_{n-m-*}(p^{-1}(y))$.
Similarly, the stalk  of the locally constant derived sheaf $\mc H^*(\mc I^{\bar p}_{fc}\mc Z^*_{E,A})$ is isomorphic to $IH^c_{n-m-*}(p^{-1}(y),p^{-1}(y)\cap A)$.
\end{proposition}
\begin{proof}
That $\mc H^*(\mc I^{\bar p}_{fc}\mc Z^*_E)$ is locally constant comes from  Corollary \ref{C: clc}. By Corollary \ref{C: H stalk}, the stalk at any point $y$ in $Y$ is $IH^c_{n-*}(E,E-p^{-1}(y))$. By excision \cite[Lemma 2.11]{GBF10}, this is isomorphic to $IH^c_{n-*}(p^{-1}(D),p^{-1}(D)-p^{-1}(y))=IH^c_{n-*}(p^{-1}(D),p^{-1}(D-y))$, where $D$ is a neighborhood of $y$ homeomorphic to $\R^m$. By \cite[Corollary 3.14]{GBF3}, there is a strong stratum-preserving fiber homotopy equivalence $p^{-1}(D)\to D\times p^{-1}(y)$. Thus $IH^c_{n-*}(p^{-1}(D),p^{-1}(D-y))\cong IH^c_{n-*}((D, D-y)\times p^{-1}(y))$. Since $(D,D-y)$ is a manifold pair, we can apply the intersection homology K\"unneth theorem (see \cite{Ki}) to get 
$\mc H^*(\mc I^{\bar p}_{fc}\mc Z^*)_y\cong IH^c_{n-m-*}(p^{-1}(y))$.

Similarly, that $\mc H^*(\mc I^{\bar p}_{fc}\mc Z^*_{E,A})$ is locally constant follows from Proposition
\ref{P: rel clc}, and by Proposition \ref{P: rel stalk},
the stalk  at $y$ is isomorphic to $IH^c_{n-*}(E,A\cup (E-p^{-1}(y)))\cong IH^c_{n-*}(p^{-1}(D), (p^{-1}(D)\cap A) \cup p^{-1}(D-y))$. 
We can then employ a relative version of \cite[Corollary 3.14]{GBF3}, using the relative lifting property, and we see that $IH^c_{n-*}(p^{-1}(D), (p^{-1}(D)\cap A) \cup p^{-1}(D-y))\cong IH^c_{n-*}(D\times p^{-1}(y),(D\times (A\cap p^{-1}(y))\cup ((D-y)\times p^{-1}(y)))\cong IH^c_{n-m-*}(p^{-1}(y),p^{-1}(y)\cap A)$. 
\end{proof}

The following corollary follows immediately from the proof of the proposition:

\begin{corollary}
Let $p: (E,A)\to Y$ be a stratified fibration pair, where $Y$ is an unfiltered $m$-manifold.  Then the locally constant derived sheaf $\mc H^*(\mc I^{\bar p}_{fc}\mc Z_E^*)$ is isomorphic to $\mc{IH}^c_{n-m-*}(p^{-1}(y))\otimes \mc O$, where $\mc O$ is the orientation sheaf on $Y$ and $\mc{IH}^c_{n-m-*}(p^{-1}(y))$ is the coefficient system with stalks $IH^c_{n-m-*}(p^{-1}(y))$ and path action determined by the stratified fibration.
Similarly, $\mc H^*(\mc I^{\bar p}_{fc}\mc Z^*_{E,A})\cong  \mc{IH}^c_{n-m-*}(p^{-1}(y),p^{-1}(y)\cap A)\otimes O$.
\end{corollary}

\begin{corollary}\label{C: manifold base PD}
Let $\mc O$ be the orientation sheaf on the $m$-manifold $Y$. Suppose that $\mc{IH}^c_{*}(p^{-1}(y))$ and $\mc{IH}^c_{*}(p^{-1}(y), p^{-1}(y)\cap A)$ are finitely generated. Then the $E_2$ terms  $H^p_c(Y;\mc{H}^q(\mc I^{\bar p}_{fc}\mc Z^*))$ and $H^p(Y;\mc{H}^q(\mc I^{\bar p}_{fc}\mc Z^*))$ of Theorem \ref{T: strat fib SS} and 
$H^p_c(Y;\mc{H}^q(\mc I^{\bar p}_{fc}\mc Z_{E,A}^*))$ and $H^p(Y;\mc{H}^q(\mc I^{\bar p}_{fc}\mc Z_{E,A}^*))$ of Theorem \ref{T: rel strat fib SS} are respectively the singular homology and cohomology groups with local coefficients 
$$H^p_c(Y; \mc{IH}^c_{n-m-q}(p^{-1}(y))\otimes \mc O)\cong 
H_{m-p}^c(Y; \mc{IH}^c_{n-m-q}(p^{-1}(y))),$$ 
$$H^{p}(Y;  \mc{IH}^c_{n-m-q}(p^{-1}(y))\otimes \mc O)\cong 
H_{m-p}^{\infty}(Y;  \mc{IH}^c_{n-m-q}(p^{-1}(y))),$$ 
$$H^p_c(Y; \mc{IH}^c_{n-m-q}(p^{-1}(y), p^{-1}(y)\cap A)\otimes \mc O)\cong 
H_{m-p}^c(Y;  \mc{IH}^c_{n-m-q}(p^{-1}(y), p^{-1}(y)\cap A)),$$ 
$$H^p(Y;  \mc{IH}^c_{n-m-q}(p^{-1}(y), p^{-1}(y)\cap A)\otimes \mc O)\cong 
H^{\infty}_{m-p}(Y;  \mc{IH}^c_{n-m-q}(p^{-1}(y), p^{-1}(y)\cap A)).$$ 
Moreover, these isomorphisms are natural.
\end{corollary}
\begin{proof}
By \cite[Theorem III.1.1]{Br}, these $E_2$ terms are just ordinary singular cohomology with local coefficients or singular cohomology with compact supports and local coefficients (or, equivalently, \v{C}ech cohomology by \cite[Corollary III.4.12]{Br}). The rest of the theorem then follows by Poincar\'e duality; see, e.g., Theorem 10.2 of \cite{Sp93} (and note that here $\bar H^*=H^*$ by the comments on \cite[p. 188]{Sp93}).
\end{proof}

In the situation of Proposition \ref{P: cocompact}, this corollary also tells us about $IH_*^{\infty}(E,A)$, when $Y$ is a manifold and $p:(E,A)\to Y$ satisfies the cocompactness property.

\vskip1cm

With manifold bases, we can also say something more about intersection homology with closed supports, provided we make some stronger assumption on the properties of the stratified fibration. Whenever we work with closed support intersection homology, we assume all spaces to be paracompact.

If we make some very strong assumptions, then we have the following:

\begin{proposition}\label{P: infty compact fibers}
Suppose that the stratified fibration $p:E\to Y$ is such that the fibers $p^{-1}(y)$ are compact  of finite homological dimension and that $Y$ is an $m$-manifold. Then $\mc{H}^q( Rp_*\mc{IS}^* )$ is locally constant with stalks $IH_{n-m-q}( p^{-1}(y))$.
\end{proposition}
\begin{proof}
We have already seen that, in general, $\mc{H}^q( Rp_*\mc{IS}^* )$
is isomorphic to $\lim_{y\in U} IH^{\infty}_{n-q}(p^{-1}(U))$. But in the present case, we can take each $U$ to be homeomorphic to $\R^m$, and then $p^{-1}(U)$ is stratum-preserving fiber homotopy equivalent to $U\times p^{-1}(y)$ by \cite[Corollary 3.14]{GBF3}. Since the fibers are assumed compact, this is a proper homotopy equivalence, and so we have, for each such $U$,  $IH^{\infty}_{n-q}(p^{-1}(U))\cong IH^{\infty}_{n-q}(U\times p^{-1}(y))$, which is isomorphic to $IH_{n-m-q}( p^{-1}(y))$ by \cite[Proposition  2.15]{GBF10}. It is clear from the methods of proof of that proposition that the local system over $U$ is constant.
\end{proof}

As for Corollary \ref{C: manifold base PD}, we have the following.

\begin{corollary}\label{C: infty compact fiber}
Under the assumptions of the preceding proposition, let $\mc O$ be the orientation sheaf on the $m$-manifold $Y$. Then the $E_2$ terms  $H^p(Y;\mc{H}^q( Rp_*\mc{IS}^* ))$
 of Theorem \ref{T: lf IH SS} are the locally finite homology groups with local coefficients $H_{m-p}^{\infty}(Y;  \mc{IH}_{n-m-q}( p^{-1}(y)))$. 
\end{corollary}

We can discard the assumption that the fibers be compact if instead we require that our stratified fibration actually be a stratified bundle. 

\begin{definition}
$p:E\to Y$ is a stratified bundle over the $m$-manifold $Y$ if for each $y\in Y$ there is a neighborhood $U\subset Y$ such that $p^{-1}(U)\cong U\times F$ for a paracompact locally-compact filtered  space $F$ of finite cohomological dimension. 
\end{definition}

\begin{proposition}
Suppose that $p:E\to Y$ is a stratified bundle over the $m$-manifold $Y$ with fibers $F$ of stratified dimension $n-m$. Then $\mc{H}^q( Rp_*\mc{IS}^* )$ is locally constant with stalks $IH^{\infty}_{n-m-q}(F)$. Furthermore, if $\mc O$ is the orientation sheaf on $Y$. Then the $E_2$ terms  $H^p(Y;\mc{H}^q( Rp_*\mc{IS}^* ))$
 of Theorem \ref{T: lf IH SS} are the locally finite homology groups with local coefficients $H_{n-p}^{\infty}(Y; \mc{IH}_{n-m-q}(X))$.
\end{proposition}
\begin{proof}
The proof is as for Proposition \ref{P: infty compact fibers} and Corollary \ref{C: infty compact fiber} using the following proposition for the local calculation.
\end{proof}

\begin{proposition}
Let $X$ be a  (not necessarily compact) paracompact filtered  space of finite cohomological dimension and stratified dimension $n$. Then $IH^{\infty}_*(X\times \R)\cong IH^{\infty}_{*-1}(X)$.
\end{proposition}
\begin{proof}
We first construct a homomorphism $h: IC_*^{c}(X)\to IC_{*+1}^{c}(X\times \R)$. We do this by choosing for each singular simplex $\sigma: \Delta\to X$ a triangulation of $\Delta\times [-1,1]$ so that $\sigma\times \text{id}:\Delta\times [-1,1]\to X\times [-1,1]$ yields a chain. These choices can be made inductively over dimension so as to make $h$ a chain map. Note that there are no problems with allowability  as $X\times [-1,1]$ is given the product filtration (see, e.g., the proof of stratum-preserving homotopy invariance of intersection homology in \cite{GBF3}).

Now, we will think of $X\times \R$ as the total space of a stratified fibration with base $X$ and fiber $\R$. The homomorphism $h$ induces a sheaf map $\mc {IS}^*_{X} \to \mc I^{\bar p}_{fc}\mc Z^{*}_{X\times \R,X\times \R^*}$, where $\R^*=\R-0$. The sheaf map is induced at the level of presheaves by the map constructed on chains. It is well-defined since chains in $IC^c_*(X-\bar U)$, $U\subset X$, map to chains in $IC^c_*((X-\bar U)\times \R)$. The hypercohomology  of  $\mc {IS}^*_{X}$ is just $IH^{\infty}_{n-*}(X)$, while that of  $\mc I^{\bar p}_{fc}\mc Z^{*}_{X\times \R,X\times \R^*}$ is $IH^{\infty}_{n+1-*}(X\times \R, X\times \R^*)$ by Proposition \ref{P: cocompact}. By \cite[Lemma 2.14]{GBF10}, $IH^{\infty}_*((X\times \R^*)_{X\times \R})=0$, so $IH^{\infty}_{n+1-*}(X\times \R, X\times \R^*)\cong IH^{\infty}_{n+1-*}(X\times \R)$. Thus to prove the proposition, it suffices to show that the sheaf map induced by $h$ is a quasi-isomorphism.

The stalk cohomology at $x$ of $\mc {IS}^*_{X}$ is $IH_{n-*}^c(X,X-x)$. This can be seen, e.g., from  Corollary \ref{C: H stalk} by treating $X$ as the trivially stratified fibration over itself. Similarly by Proposition \ref{P: rel stalk}, the stalk at $x$ of $\mc I^{\bar p}_{fc}\mc Z^{*}_{X\times \R,X\times \R^*}$ is $IH^c_{n+1-*}(X\times \R,X\times \R-x\times 0)=IH^c_{n+1-*}((X,X-x)\times (\R,\R^*))$. By the intersection homology K\"unneth theorem \cite{Ki}, this is $IH^c_{n-*}(X,X-x)$. It only remains to see that the isomorphism is indeed induced by $h$. But this follows as for the ordinary K\"unneth theorem since the oriented simplex $[-1,1]$ generates $H_1(\R,\R-0)$, and our $h$ can be interpreted as the cross product $\xi\to \xi\times [-1,1]$. One need only extend the K\"unneth theorem in \cite{Ki} to the relative cases. (Alternatively, we have the long exact sequence 
\begin{equation*}
\begin{CD}
@>>>&IH^c_*((X,X-x)\times \R^*) &@>>>&IH^c_*((X,X-x)\times \R)&@>>>&IH^c_*(X\times \R, X\times \R-x\times 0)&@>\bd_*>>,
\end{CD}
\end{equation*}
where the last term is correct since the excisive couple property of open sets holds in intersection homology by \cite{GBF10}. By stratum-preserving homotopy equivalence of compactly supported intersection homology, this is isomorphic to
\begin{equation*}
\begin{CD}
@>>>&IH^c_*(X,X-x)\oplus IH^c_*(X,X-x) &@>>>&IH^c_*(X,X-x)&@>>>&IH^c_*(X\times \R, X\times \R-x\times 0)&@>>>.
\end{CD}
\end{equation*}
It follows that $IH^c_*(X\times \R, X\times \R-x\times 0)\overset{\cong}{\to}IH^c_{*-1}(X,X-x)$ and that the map induced by $h$ is a right inverse to $\pi\bd_*$ by construction of $h$, where $\pi$ is projection to the $IH^c_*((X,X-X)\times (0,\infty)$  summand.) 
\end{proof}

\begin{corollary}
With the assumptions of the Proposition, $IH^{\infty}_*(X\times \R^k)\cong IH^{\infty}_{*-k}(X)$.
\end{corollary}
\begin{proof}
By induction.
\end{proof}

\section{Neighborhoods}\label{S: neighborhoods}

We can now employ our results concerning the intersection homology of stratified fibrations to say something about the intersection homology of neighborhoods or pure subsets of stratified spaces. At this point, the reader may want to review the definitions in Section \ref{S: def MHSSs}.

A \emph{pure subset} $Y$ of a filtered space $X$ is a closed subspace that is the union of connected components of strata. In other words, $Y$ should contain each component of a stratum that it intersects. Examples include the skeleta of $X$. We will mostly limit ourselves to pure subsets of manifold homotopically stratified spaces (MHSSs). 

We wish to study the intersection homology of neighborhoods of pure subsets of MHSSs. The most general types of neighborhoods we shall consider are nearly stratum-preserving deformation retract neighborhoods (NSDRNs). We recall from Section \ref{S: def MHSSs}
that an NSDRN of $Y$ is a neighborhood $N$ of $Y$ such that there exists a strong deformation retraction $r: N\times I\to N$ that retracts $N$ to $Y$ and such that $r$ is stratum-preserving on $N\times [0,1)$. Examples include 
\begin{itemize}
\item stratified tubular neighborhoods, retracted in the obvious way, 
\item stratified PL regular neighborhoods, retracted along join lines (see \cite[\S 70]{MK}),
\item stratified mapping cylinders (see \cite{Hu99a, CS95}), retracted down the cylinder,
\item Hughes's approximate tubular neighborhoods \cite{Hu02}.
\end{itemize}
The last item generalizes the others and is the most natural candidate for the role of ``regular'' neighborhoods in the category of MHSSs. See \cite{Hu02} for more about these neighborhoods, including an existence proof (subject to the condition that all non-minimal strata have dimension $\geq 5$). 

We will prove that approximate tubular neighborhoods are NSDRNs and that they are \emph{outwardly stratified tame}. These proofs concerning properties of approximate tubular neighborhoods are collected separately in Section \ref{S: ATN}.

\begin{lemma}\label{L: ATN=NSDRN}
Suppose that $N$ is an approximate tubular neighborhood of $Y$ in the MHSS $X$. Then there is a nearly stratum-preserving deformation retraction taking $N$ into $Y$ rel $Y$.
\end{lemma}
\begin{proof}
See Section \ref{S: ATN}
\end{proof}

Thus we have the following immediate corollary by the definition of NSDRNs.
\begin{corollary}\label{C: ATN=NSDRN}
Approximate tubular neighborhoods in MHSSs are NSDRNs. 
\end{corollary}

\subsection{$IH^c_*$ of NSDRNs}

Now, to study the intersection homology of NSDRNs, we use the fact that NSDRNs are stratum-preserving homotopy equivalent to mapping cylinders of stratified fibrations. In fact, $N$ 
is stratum-preserving homotopy equivalent rel $Y$ to the mapping cylinder $M$ of the holink evaluation $\pi: \hl_s(N,Y)\to Y$. This proposition is essentially found in the work of Chapman \cite{Ch79} and Quinn \cite{Q1}; a detailed (and corrected) proof is given in the appendix to \cite{GBF3}. The map $\pi$ is a stratified fibration by \cite{Hug}.

Also, it is shown in \cite[Proposition 3.3]{GBF5} that the mapping cylinder collapse of a mapping cylinder of a stratified fibration is itself a stratified fibration. Thus, $N$ is stratum-preserving homotopy equivalent rel $Y$ to the stratified fibration $M\to Y$. Furthermore, since $N$ is a neighborhood of $Y$ in the MHSS $X$, $\hl_s(X,Y)\sim_{sphe}\hl_s(N,Y)$, since each is stratum-preserving homotopy equivalent to the holink of small paths, $\hl_s^{\delta}(X,Y)$, for some sufficiently chosen $\delta$ (see \cite{Q1}). Let $\mf M$ be the mapping cylinder of the holink evaluation $\hl_s(X,Y)\to Y$.

It is not hard to see now that we have stratum-preserving homotopy equivalences of pairs 
$$(N,N-Y)\sim_{sphe} (M,M-Y)\sim_{sphe}(M,\hl_s(N,Y))\sim_{sphe}(\mf M, \hl_s(X,Y))\sim_{sphe}(\mf M,\mf M-Y).$$

By Proposition 2.1 of \cite{GBF3}, $IH^c_*$ is a stratum-preserving homotopy invariant. Thus we have the following:

\begin{proposition}
Let $Y$ be a pure subset of the MHSS $X$, and let $N$ be an NSDRN  of $Y$. Let $M$ and $\mf M$ be the mapping cylinders of the evaluations $\hl_s(N,Y)\to Y$ and $\hl_s(X,Y)\to Y$. Then the long exact intersection homology sequences with compact supports of the following pairs of spaces are isomorphic: \begin{itemize}\item $(N,N-Y)$,
\item $(M,M-Y)$, 
\item $(\mf M,\mf M-Y)$, 
\item $(M,\hl_s(N,Y))$,  and  
\item $(\mf M, \hl_s(X,Y))$.
\end{itemize} 
Furthermore, the pairs $(M,M-Y)$ and $(\mf M,\mf M-Y)$ are stratified fibration pairs over $Y$ in the sense of Definition \ref{D: strat fib pair}.
\end{proposition}
\begin{proof}
The proof follows from the stratum-preserving homotopy equivalences listed above. The only part needing comment is to note that $(M,M-Y)\to Y$ and $(\mf M,\mf M-Y)$ are stratified fibrations \emph{as pairs} because $Y$ is a union of components of strata and $M\to Y$ and $\mf M\to Y$ are already stratified fibrations.
\end{proof}

The results of the preceding sections apply to the stratified fibration pairs $p:(M,M-Y)\to Y$ and $p:(\mf M, \mf M-Y) \to Y$ and demonstrate that the compactly supported intersection homology groups of $M$, $M-Y$, and $(M, M-Y)$, which are isomorphic to the intersection homology groups of $N$, $N-Y$ and $(N,N-Y)$, are the abutments of spectral sequences whose $E^2$ terms are the compactly supported cohomology of stratified cohomologically locally constant sheaves on $Y$:

\begin{corollary}\label{C: NSDRN IHC}
$IH^c_*(N)$, $IH^c_*(N-Y)$, and $IH^c_*(N,N-Y)$ are the abutments of spectral sequences with $E^2$ terms $H^p_c(Y;\mc{H}^q(\mc I^{\bar p}_{fc}\mc Z_M^*))$, $H^p_c(Y;\mc{H}^q(\mc I^{\bar p}_{fc}\mc Z_{M-Y}^*))$, and $H^p_c(Y;\mc{H}^q(\mc I^{\bar p}_{fc}\mc Z_{M,M-Y}^*))$.
\end{corollary}
\begin{proof}
$N$ and $N-Y$ are stratum preserving homotopy equivalent to stratified fibrations, and compact support intersection homology is an invariant with respect to such equivalences. We then apply the results of Section \ref{S: strat fib} to these stratified fibrations, noting that all spaces are metric and $Y$ is locally-compact of finite cohomological dimension (as a closed subset of $X$).
\end{proof}

\subsection{$IH^{\infty}_*$ of NSDRNs}\label{S: N closed}

We can also use the tools of the earlier sections to say something about closed support intersection homology of neighborhoods, though, as we have seen before, we cannot take a direct route since the stratum-preserving homotopy equivalence from $N$ to the mapping cylinder $M$ of  $\hl_s(N,Y)\to Y$ is not proper. Nonetheless, we can obtain information about $IH^{\infty}_*(N)$ by showing that it is isomorphic to $IH^{\infty}_*(N,N-Y)$ under certain assumptions on the NSDRN $N$. These assumptions will hold if $N$ is an approximate tubular neighborhoods (and so a regular neighborhood, tubular neighborhood, or mapping cylinder neighborhood). From there, we show that  $IH^{\infty}_*(N,N-Y)\cong IH^{fc}_*(M, M-Y)$, which we have already treated above. 

We first set forth the extra property we will need on our NSDRNs. This is a certain tameness condition that says, essentially, that sets in $N-Y$ that do not intersect $Y$ can be pushed away off the far end of $N$.

\begin{definition}\label{D: ost}
We say that the neighborhood $N$ of $Y$ is \emph{outwardly stratified tame} if it possesses the following property: For any metric space $Z$ and proper map $f: Z\to N$ such that $f(Z)\subset N-Y$, there exists a proper stratum-preserving homotopy $H: Z\times [0,\infty)\to N$ such that $H(Z\times [0,\infty))\subset N-Y$ and $H|_{Z\times 0}=f$.
\end{definition}

Approximate tubular neighborhoods in MHSSs are outwardly stratified tame:

\begin{proposition}\label{P: ATN=OST}
Approximate tubular neighborhoods of pure subsets of MHSSs are outwardly stratified tame.
\end{proposition}
\begin{proof}
See Section \ref{S: ATN}, below.
\end{proof}

We see why this property is useful by showing that if $N$ is outwardly stratified tame then $IH^{\infty}_*(N)\cong IH^{\infty}_*(N,N-Y)$. The group $IH^{\infty}_*(N,N-Y)$ is usually easier to work with then $IH^{\infty}_*(N)$ as we have seen above in Proposition \ref{P: cocompact}, in \cite{GBF10}, and as we will demonstrate further below.

\begin{proposition}
Let $N$ be an outwardly stratified tame NSDRN of $Y$. Then $IH^{\infty}_*(N)\cong IH^{\infty}_*(N,N-Y)$. 
\end{proposition}
\begin{proof}
Recall that $IH^{\infty}_*(N,N-Y)$ is the homology of the quotient chain complex $IC_*^{\infty}(N)/IC_*^{\infty}((N-Y)_N)$ and that $IC_*^{\infty}((N-Y)_N)$ consists of allowable chains with supports in $N-Y$ that are locally finite  in $N$.  We will show that $IH_*^{\infty}((N-Y)_N)=0$, which will prove the proposition by the long exact sequence of the pair.

Let $\xi$ be a cycle in $IC_i^{\infty}((N-Y)_N)$. Each simplex $\sigma_j$ in $\xi$ is a map $\Delta^i_j\to N-Y$, where $\Delta^i_j$ is a copy of the standard model simplex. So we can think of the chain $\xi$ as being determined by a map $f: \amalg_j \Delta_j^i\to N-Y$, along with the coefficient information. (Note that the index set for $j$ need not be countable, and we can treat $\amalg_j \Delta_j$ as a metric space by setting the distance between connected components to be infinite.) The  map $f$ is proper because any compact set in $N$ can intersect only a finite number of simplices by the local finiteness condition on $\xi$. 

Let $|\xi|$ be the support of $\xi$, and let $g:|\xi|\to N$ denote the inclusion. We may identify $|\xi|$ with the image of $f$ by assuming that we have eliminated from $\xi$ all singular simplices whose coefficients add to $0$ (this presents no conflict with allowability). Let $f'$ denote $f$ with codomain restricted to $|\xi|$ so that $f=gf'$. $f'$ is clearly proper since $f$ is. $g$ is also proper, since for any compact set $C$ in $N$, $C$ intersections only finitely many simplices of $\xi$ and thus $g^{-1}(C)=|\xi|\cap C$ is the intersection of $C$ with the image of finitely many simplices under $f$. 

By the outward stratified tameness of $N$, there is an extension of $g$ to a stratum-preserving homotopy $H: |\xi|\times [0,\infty)\to N-Y$ that is proper as a map to $N$. We can use this homotopy to build a homology from $\xi$ to the empty chain. In particular, we triangulate $\amalg_j \Delta_j^i\times [0,\infty)$ compatibly so that the products with $[0,\infty)$ of common faces of singular simplices receive the same triangulations, and then for each simplex $\Delta_j^i$, we  consider the composition $H\circ (f'\times \text{id}_{[0,\infty)}):\Delta_j^i\times [0,\infty)\to N-Y$. Together with
the original coefficient information, these maps are used to build a new chain $\Xi$ such $|\Xi|\subset N-Y$ and $\bd \Xi=\xi$. Note that any agreeing faces of simplices are homotoped together and that the chosen compatible triangulations insure that $\xi$ does not acquire any new unwanted boundaries as it is pushed to infinity.
Technical details of an analogous construction are given in the proofs of \cite[Proposition 2.7 and Lemma 2.14]{GBF10}. 

The allowability of $\Xi$ follows as for any homology built from a stratum-preserving homotopy (see \cite{GBF3}). To check the locally finiteness, we note that if $x\in N$, then since $N$ is locally-compact, $x$ has a neighborhood $K$ that is compact. But then since $H$ and $f'\times \text{id}_{[0,\infty)}$ are proper, $(H\circ (f'\times \text{id}_{[0,\infty)}))^{-1}(K)$ can intersect only finitely many of the polygonal $i$-simplices in the triangulation of $\amalg_j \Delta^i_j\times [0,\infty)$. Thus $\Xi$ is locally-finite.

\end{proof}

Roughly speaking, the preceding  proposition says that we can take cycles in $N-Y$ and ``push them off to infinity''.

We have seen that if $N$ is outwardly stratified tame, then $IH^{\infty}_*(N)$ reduces to $IH^{\infty}_*(N,N-Y)$. 
Now we show that $IH_*^{\infty}(N, N-Y)\cong IH_*^{fc}(M,M-Y)$, which we have already studied as $(M, M-Y)$ is a stratified fibration pair.

\begin{theorem}\label{T: transition}
If $N$ is an NSDRN of the pure subset $Y$ in an MHSS and  $M$ is the mapping cylinder of the holink evaluation $\pi: \hl_s(N,Y)\to Y$, then $IH^{\infty}_*(N,N-Y)\cong IH_*^{fc}(M,M-Y)$.
\end{theorem}

The proof requires several lemmas.

\begin{lemma}
Let $\delta: Y\to (0,\infty)$ be a continuous function. Let $M^{\delta}$ be the mapping cylinder of $\pi: \hl^{\delta}_s(N,Y)\to Y$. Then the inclusion $M^{\delta}\into M$ induces an isomorphism $IH_*^{fc}(M,M-Y)\cong IH_*^{fc}(M^{\delta}, M^{\delta}-Y)$. If $\delta_1,\delta_2$ are two such functions with $\delta_2\leq \delta_1$, then inclusion induces an isomorphism $IH_*^{fc}(M^{\delta_2},M^{\delta_2}-Y)\cong IH_*^{fc}(M^{\delta_1}, M^{\delta_1}-Y)$.
\end{lemma}
Note: although $p: M^{\delta}\to Y$ may not be a stratified fibration, the definition of the support condition on $IH_*^{fc}(M^{\delta}, M^{\delta}-Y)$  generalizes in the obvious way.
\begin{proof}
By \cite[Lemma 2.4]{Q1} the inclusion $\hl^{\delta}_s(N,Y)\into \hl_s(N,Y)$ is a stratum-preserving homotopy equivalence, obtained by shrinking paths.  It follows then that $M^{\delta}\into M$ is a stratified fiberwise homotopy equivalence and similarly for the pairs $(M,M-Y)$ and $(M^{\delta}, M^{\delta}-Y)$. 

These equivalences can then be used to construct the desired intersection homology isomorphism provided we are careful since we are using $fc$ intersection homology and not compactly supported intersection homology. However, the support condition is equivalent to asking that every point  $y\in Y$ has a neighborhood $U$ such that $p^{-1}(U)$ intersects only a finite number of simplices of any cycle representing a given intersection homology class. But this property is clearly preserved under the homotopy equivalence since all maps and homotopies are fiber-preserving. 

The second statement follows similarly.
\end{proof}

We now define a specific $\delta_1$ as follows: for each point $y\in Y$, assign an open neighborhood $U_y$ of $y$ in $N$ such that $\bar U_y$ is compact. Let $\mc V$ be a locally-finite refinement of $\{U_y\}\cup (N-Y)$. This is possible since $N$ is paracompact and locally compact. For each $y\in Y$, let $V_{(y)}$ be a set of $\mc V$ containing $y$. Now use the metric on $y$ to define $\delta_1$ as a continuous function such that for each $y$ in $Y$, $B_{\delta_1(y)}(y)\subset V_{(y)}$, where $B_{\delta_1(y)}(y)$ is the closed ball of radius $\delta_1(y)$ about $y$ in $N$. Note that each $V_{(y)}$ must be contained in some $U_y$, and so $\bar V_{(y)}$ is compact.

We also build a certain neighborhood of $Y$ in $N$. Let $R$ be a nearly stratum-preserving deformation retraction of $N$ to $Y$, and let $R_x$ be the retraction path of $x$.  Let $Q'=\{x\in N\mid R_x\subset B_{\delta_1(R_x(0))}(R_x(0))\}$. In other words, $Q'$ is the set of points whose retraction paths lie within a distance $\delta_1(R_x(0))$ of their endpoints $R_x(0)$ in $Y$. Note:
\begin{enumerate}
\item $Y\subset Q'$ since for each $y\in Y$, $R_y$ is constant and holds $y$ fixed.

\item $Q'$ is open. This follows easily from the continuity of $R$ and $\delta_1$.  

\item If $x\in Q'$, then $R_{x}\subset V_{(R_x(0))}$ from the definitions.
\end{enumerate}
Let $Q=Q'\cap (\cup_y V_{(y)})$, where the $V_{(y)}$ are as in the definition of $\delta_1$.

In what follows, we will also use some of the notation of Proposition A.1 of \cite{GBF3}, which demonstrates the stratum-preserving homotopy equivalence of $M$ and $N$ (note, in \cite{GBF3} $M$ is referred to as $I_{\pi}$ and our $Y$ is there called $X$).    We let $f:N\to M$ and $g:M\to N$ be the homotopy inverses. They are defined by $f(y)=y$ and $g(y)=y$ for $y\in Y$. For $x\in N$, $f(x)=(R_x,d(x,Y))$, where the distance $d$ has been suitably rescaled so that $d(x,Y)<1$ for all $x\in N$. The pair $(R_x, d(x,Y))$ represents coordinates in the mapping cylinder of the holink evaluation $\pi:\hl_s(N,Y)\to Y$. The map $g$ is more complicated; it takes a point $(\omega,t)$ in the mapping cylinder to a certain point along the path $\omega$. See \cite{GBF3} for a complete description. 

Throughout the following, we will use the same symbols $f$ and $g$ to denote both the maps and how they act on chains. We avoid giving them special names like $g_*$ or $f_*$ since often the ``chain maps'' are only well-defined for the specific chains we will consider below.

\begin{lemma}\label{L: g map}
For any $\delta\leq \delta_1$, $g$ induces a homomorphism $g:IH^{fc}_*(M^{\delta},M^{\delta}-Y)\to IH_*^{\infty}(N, N-Y)$
\end{lemma}
\begin{proof}
Recall that the map $g$, as defined in Proposition A.1 of \cite{GBF3} takes $Y$ identically to itself and takes each $(\omega, s)\in M-Y$ to a point in the path $\omega$. Since $g$ is stratum-preserving, the only concern is that images of cycles be locally-finite.

So let $\xi$ be a cycle in $IH_*^{fc}(M^{\delta},M^{\delta}-Y)$. By the definition of $\delta_1$, it is clear that $g(|\xi|)\subset \cup V_{(y)}$. So let $x\in \cup V_{(y)}$, and suppose in particular that $x$ is in the specific element $V$ of $\{V_{(y)}\}$. A simplex $\sigma$ of $\xi$ can have $g(|\sigma|)\cap V\neq \emptyset$ only if there is a point $z\in p(|\sigma|)\subset Y$ such that $V_{(z)} \cap V\neq \emptyset$ since, by definition of $\delta_1$, for each $(\omega,s)\in M^{\delta}$, $\omega(I)\subset V_{\pi(\omega)}$. But the number of such $V_{(x)}$ (including $V$ itself) is finite, since $\mc V$ is locally-finite and $\bar V$ is compact. Call this collection $V_1,\ldots, V_k$. But $W=\cup_{i=1}^k \bar V_i=\overline{\cup_{i=1}^k V_i}$ is compact and hence so is $W\cap Y$. But this implies that the number of simplices in $\xi$ intersecting $p^{-1}(W \cap Y)$ is finite by using the fc support condition.  Since these are the only simplices whose image under $g$ can intersect $V$, $V$ is a neighborhood of $x$ intersecting only a finite number of simplices of $g(\xi)$. 
\end{proof}

Let $IH^{\infty}_*(Q_N, (Q-Y)_N)=H_*(IC^{\infty}_*(Q_N)/ IC^{\infty}_*((Q-Y)_N))$. Recall that this means that that all chains must be locally finite in $N$ (not just in $Q$). 

\begin{lemma}\label{L: Q iso}
The map $IH^{\infty}_*(Q_N, (Q-Y)_N)\to IH^{\infty}_*(N,N-Y)$ induced by inclusion is an isomorphism.
\end{lemma}
\begin{proof}
Let $[\xi]\in IH^{\infty}_i(N,N-Y)$. We claim that $[\xi]=[\eta]$ for some $\eta\in IC_i^{\infty}(N)$ with $|eta|\subset Q$. This is essentially a subdivision and excision argument as in the proof of Proposition \ref{P: cocompact}. First we find a subdivision of $\xi'$ of $\xi$ such that each simplex intersecting $Y$ lies in $Q$. This can be done inductively by successively subdividing skeleta rel the lower skeleta. Then we fix the (not necessarily allowable) chain $\Xi$ consisting of simplices of $\xi'$ (with their coefficients) that intersect $Y$. Then we further subdivide $\xi'$ rel $\Xi$ until all simplices sharing vertices with $\Xi$ also lie in $Q$. Then we let $\eta$ equal $\Xi$ plus all simplices (with coefficients) of the new subdivision that share vertices with $\Xi$. Then $\eta$ is allowable, $|\eta|\subset Q$, and $[\eta]=[\xi]$. Similarly, if $[\xi]\in IH^{\infty}_*(Q_N, (Q-Y)_N)$ is homologous to $0$ in $IH^{\infty}_*(N,N-Y)$, then we can perform a similar subdivision and excision (excising away from $|\xi|\cup Y$) on the null-homology, either rel $\xi$ or letting $\xi$ get subdivided and noting that the subdivided chain represents the same element of $IH^{\infty}_*(Q_N, (Q-Y)_N)$. The details are similar to (though simpler than) those of Proposition \ref{P: cocompact}, and more detailed versions of similar arguments can be found in \cite{GBF10}. 
\end{proof}

\begin{lemma}\label{L: f map}
$f$ induces a map $IH^{\infty}_*(Q_N,(Q-Y)_N)\to IH_*^{fc}(M^{\delta_1},M^{\delta_1}-Y)$. 
\end{lemma}
\begin{proof}
We show that $f$ (as a chain map), takes a cycle $\xi\in IC^{\infty}_i(Q_N,(Q-Y)_N)$ into $IC_i^{fc}(M^{\delta_1},M^{\delta_1}-Y)$. First, by the definition of $Q$, $f$ takes all of $Q$ into $M^{\delta_1}$, and the fact that $f$ is stratum preserving guarantees the allowability of $f(\xi)$. The only thing then to check is the support condition. So let $y\in Y$, and consider $Z=V_{(y)}\cap Y$. The simplices $\sigma$ of $\xi$ such that  $f(|\sigma|)$ intersect $p^{-1}(Z)$ are those containing points $x$ such that $R_x(0)\in Z$. By definition of $\delta_1$, for each $z\in Z$, $B_{\delta_1(z)}(z)\subset V_{(z)}$, so for such $x$, $R_x(I)\subset V_{(z)}$. Furthermore, since the collection of sets $\mc V$ is locally finite and each $V_{(z)}$ has compact closure, there are only a finite number of $V_{(z)}$ that intersect $Z$, and these must contain all paths in $\cup_{z\in Z} B_{\delta_1(z)}(z)$. So each $\sigma$ in $\xi$ such that $pf(|\sigma|)\cap Z\neq \emptyset$ must intersect the compact set $\cup_{z\in Z} \bar V_{(z)}$. Thus there are a finite number of such simplices since $\xi$ is locally-finite. This implies the support condition.
\end{proof}

\begin{lemma}
The composition $gf$ of the maps $f$ of Lemma \ref{L: f map} and $g$ of Lemma \ref{L: g map} induce a map $IH_*^{\infty}(Q_N, (Q-Y)_N)\to IH_*^{\infty}(N, N-Y)$ that is homologous to the inclusion. In particular, it is an isomorphism. 
\end{lemma}
\begin{proof}
The second statement follows from the first and Lemma \ref{L: Q iso}.

Let $[\xi]\in IH_*^{\infty}(Q_N, (Q-Y)_N)$. Let $\xi$ be a representative of $[\xi]$. By Lemma \ref{L: f map}, $f$ takes $\xi$ to a cycle in $IC^{fc}_*(M^{\delta_1},M^{\delta_1}-Y)$, and by  Lemma \ref{L: g map}, $g$ induces a homomorphism $g:IH^{fc}_*(M^{\delta_1},M^{\delta_1}-Y)\to IH_*^{\infty}(N, N-Y)$. So the composition $gf$ takes $\xi$  to an allowable cycle in $IH_*^{\infty}(N, N-Y)$. We claim that this cycle is intersection homologous to $\xi$ in $IH^{\infty}_*(N,N-Y)$. The homology can be constructed from the stratum-preserving homotopy $gf\sim_{s.p.} \text{id}$ (as in \cite{GBF3} or \cite{GBF10}) provided we show that no local-finiteness conditions will be violated. Let's call the homotopy $H: N\times I\to N$. As constructed in \cite{GBF3}, $H$ simply retracts the paths $R_x$ along themselves in an appropriate manner.

So first suppose $z \in \bar V_{(y)}$ for some $y\in Y$ and that $\sigma$ is a simplex of $\xi$. Then $H(|\sigma|\times I)\cap V_{(y)}\neq \emptyset$ only if there is some point $x\in |\sigma|$ such that $R_x\cap V_{(y)}\neq \emptyset$. But since $x\in Q$, we know that $R_x\subset V_{(R_x(0))}$, and since each $V\in \mc V$ such that $V\cap Y\neq \emptyset$ has compact closure, there can be only a finite number of such $V$ that intersect $V_{(y)}$. Also because of the compact closures, the union of these $V$s intersects only a finite number of simplices of $\xi$. Thus only a finite number of simplices in the homology from $\xi$ to $gf(\xi)$ built from the homotopy $H$ can intersect $V_{(y)}$, a neighborhood of $z$. 

Finally suppose $z\in N$, $z\notin  V_{(y)}$ for all $y\in Y$. Then $z$ is still in some  $V\in \mc V$, and since $\mc V$ is locally finite, $V$ intersects only a finite number of $V_{(y)}$. But we know that $H(|\xi|\times I)\subset \cup_{y\in Y} V_{(y)}$ and from the last paragraph that each $\bar V_{(y)}$  intersects only finitely many simplices of the homology built from $H$. Thus $V$ is a neighborhood of $z$ intersecting only finitely many such simplices. Thus the chain providing the homology is locally finite.
\end{proof}

Let $\delta_2:Y\to (0,\infty)$ be such that $\delta_2<\delta_1$ and $B_{\delta_2(y)}(y)\subset Q$. 

\begin{lemma}
The composition $IH^{fc}_*(M^{\delta_2}, M^{\delta_2}-Y)\to IH^{\infty}_i(Q_N,(Q-Y)_N) \to IH^{fc}_*(M^{\delta_1}, M^{\delta_1}-Y)$ induced by $f$ and $g$ is the same as that induced by inclusion, and hence it is an isomorphism.
\end{lemma}
\begin{proof}
$g$ induces a homomorphism $IH^{fc}_*(M^{\delta_2}, M^{\delta_2}-Y)\to IH^{\infty}_i(N,N-Y)$ by Lemma \ref{L: g map}, but it is clear from the definition of $\delta_2$ that any cycle representing an element of the former group has its support mapped into $Q$. The map induced by $f$ is that of Lemma \ref{L: f map}. So if $\xi$ is a cycle in $IC^{fc}_*(M^{\delta_2}, M^{\delta_2}-Y)$, $fg(\xi)$ is a cycle in $IC^{fc}_*(M^{\delta_1}, M^{\delta_1}-Y)$. As in the preceding lemma, we can use the fact that $fg: M\to M$ is stratum-preserving homotopic to the identity to construct a homology between $\xi$ and $fg(\xi)$ in $IH^{fc}_*(M^{\delta_1}, M^{\delta_1}-Y)$ provided we can show that we do not violate the support condition. All other allowability concerns are taken care of since the homotopy is stratum-preserving. We do not even need to worry if the homotopy leaves $M^{\delta_1}$ since we know the inclusion $(M^{\delta_1},M^{\delta_1}-Y)\into (M, M-Y)$ induces an isomorphism on $IH^{fc}_*$, so we might as well work in $(M, M-Y)$ at this point. 

Now, let $H: M\times I\to M$ denote the homotopy from $fg$ to the identity constructed in \cite{GBF3}. This homotopy is fixed on $Y$, but it is rather complicated on $M-Y$. We do note the following, however: given $(\omega, s)\in M-Y$, the projection of $H((\omega,s), u)$ to $\hl_s(N,Y)$ is a path that lies in $\cup_{v\in [0,1]} R_{\omega(v)}(I)$. In other words, we get a path that lives in the collection of retraction paths under $R$ of points in the path $\omega$. This will be the important thing for our proof. Given a path  $\omega\in \hl_s(N,Y)$, let $T(\omega)=\cup_{v\in[0,1]}R_{\omega(v)}(0)$. This set is compact since it is the image of $[0,1]$ under the composition of continuous maps $\omega$ and $R(\cdot,0)$.

So let $y\in Y$, let $Z=V_{(y)}\cap Y$, and consider $p^{-1}(Z)\subset M$. We ask for which simplices $\sigma$ of $\xi$ will $H(|\sigma|\times I)$ intersect $p^{-1}(Z)$. 
If $x$ is a point in $\sigma$, then $x=(\omega,s)$ for some path $\omega$ and $H(x\times I)\cap p^{-1}(Z)\neq \emptyset$ only if $T(\omega)\cap Z\neq \emptyset$. Since $x\in M^{\delta_2}$, $\omega([0,1])\subset V_{(\omega(0))}$. If $T(\omega)\cap Z\neq \emptyset$, then $R_{\omega(t)}(0)\in Z$ for some $t\in [0,1]$. But then $R_{\omega(t)}\subset V_{(R_{\omega(t)}(0))}$, $V_{(R_{\omega(t)}(0))}\cap Z\neq \emptyset$, and $V_{(R_{\omega(t)}(0))}\cap V_{(\omega(0))}\neq \emptyset$, again since for any $q
\in Q$, $R_q(I)\subset V_{(q(0))}$. So $H(x\times I)\cap p^{-1}(Z)\neq \emptyset$ 
only if $p(x)$ is in a set $V\in \{V_{(z)}\mid z\in Y\}$ that intersects another such set that intersects $Z$. Since the collection $\{V_{(z)}\mid z\in Y\}$ is locally-finite and each element has compact closure  in $N$, the subcollection of elements that intersect elements that intersect $Z$, say $\mc W$, is finite and $K=Y\cap (\cup_{V\in \mc W} \bar V)$ is compact. Thus $p^{-1}(K)$ intersects only a finite number of simplices of $\xi$, by the support condition, and only the finite number of simplices in the homology induced by these can intersect $p^{-1}(Z)$. Thus the homology is allowable in $IH^{fc}_*(M,M-Y)$. 
\end{proof}

Putting these lemmas together completes the proof of Theorem \ref{T: transition}. \hfill\qedsymbol

Thus composing the results of this Section with Theorem \ref{T: rel strat fib SS}, we haven proven the following:

\begin{theorem}\label{T: closed SS}
If $N$ is an outwardly stratified tame NSDRN of the pure subset $Y$ of the MHSS $X$, then $IH^{\infty}_{n-*}(N)\cong IH^{fc}_{n-*}(M, M-Y)$, where $M$ is the mapping cylinder of the evaluation $\hl_s(N,Y)\to Y$. This group is the abutment of a spectral sequence with $E^2$ terms $H^p(Y;\mc{H}^q(\mc I^{\bar p}_{fc}\mc Z_{M,M-Y}^*))$. 
\end{theorem}

\section{Calculations in pseudomanifolds}\label{S: pms}

In the preceding section we saw how to apply the spectral sequences associated to stratified fibrations to the question of computing intersection homology of neighborhoods, specifically NSDRNs and approximate tubular neighborhoods of pure subsets of MHSSs.  The $E^2$ terms of these spectral sequences were expressed in terms of the cohomology of $Y$ with coefficients in stratified cohomologically locally constant sheaves whose stalks were given by the intersection homology groups in the path spaces $\hl_s(N,Y)$. 

In this section, we limit ourselves to pseudomanifolds, which have more concrete local structures than MHSS's, namely distinguished neighborhoods.

We recall that an $n$ dimensional stratified pseudomanifold is an $n$ dimensional MHSS such that $X-X^{n-1}$ is dense in $X$ and such that for each point $x\in X_m$, there is a \emph{distinguished neighborhood} stratified homeomorphic to $\R^m\times cL$, where the \emph{link} $L$ is an $n-m-1$ dimensional stratified pseudomanifold. These are the spaces on which intersection homology was first defined by Goresky and MacPherson \cite{GM1, GM2}, though with the extra stipulation that $X^{n-1}=X^{n-2}$.

We  wish to reinterpret the stalks of our sheaves on $Y$ in terms of the intersection homology of the actual geometric links, which are much more accessible than the homotopy links. In the topological category, links are not uniquely defined, but their intersection homology is and, either way, our computations will work for any allowable link.

Let $\pi:\hl_s(N,Y)\to Y$ be the holink evaluation and let $p:M\to Y$ be the mapping cylinder projection induced by $\pi$. Suppose $y\in Y$, and fix a distinguished neighborhood $U=\R^m\times cL$ of $y$. Let $k=\dim(L)=n-m-1$. Based on our above work, we want to compute $IH^c_*(M, M-p^{-1}(y))$, $IH^c_*(M-Y,M-(Y\cup p^{-1}(y)))$, and $IH^c_*(M, M-y)$. 

The simplest of these computations is the last one. We let $v$ stand for the vertex of $cL$.

\begin{proposition}\label{P: easy pm}
$$
IH^c_*(M, M-y)\cong IH^c_*(N, N-y)\cong IH^c_{*-m}(cL, cL-v).$$
If $\bar p$ is a traditional perversity and the coefficient system is constant, this is isomorphic to
\begin{equation*}
\begin{cases}
0,              &    *<n-\bar p(n),\\
IH^c_{*-1}(L),  &    *\geq n-\bar p(n).
\end{cases}
\end{equation*}
\end{proposition}
\begin{proof}
The second isomorphism is shown in the proof of \cite[Proposition 2.20]{GBF10}, and the last statement is well-known (see \cite{Ki}). So we focus on the first.

But we already know that the stratum-preserving homotopy equivalence $(M,M-Y)\sim_{sphe} (N,N-Y)$ is fixed on $Y$, so clearly $(M,M-y)\sim_{sphe} (N,N-y)$. 
\end{proof}

Next we prove the following theorem, which describes $IH^c_*(M, M-p^{-1}(y))$ and $IH^c_*(M-Y,M-(Y\cup p^{-1}(y)))$ in terms of the intersection homology of the link $L$ and the intersection homology of an NSDRN $N^T$ of $L^T=L\cap Y$ in $L$ ($T$ for ``tangential''). If $X$ is a PL stratified pseudomanifold, then so too will be $L$ and the neighborhood $N^T$ is guaranteed to exist (just take a regular neighborhood of $L^T$). If $X$ is a topological pseudomanifold, we will need to assume the existence of $N^T$. Of course $N^T$ will also be guaranteed to exist if the dimension conditions of Hughes's Approximate Tubular Neighborhood Theorem apply.

\begin{theorem}\label{T: pm stalk}
Let $Y$ be a pure subset of the $n$-dimensional topological stratified pseudomanifold $X$. Let $y\in Y$ such that $y\in X_m=X^m-X^{m-1}$, $y$ has a neighborhood $\R^m\times cL$, and $\dim(L)=k=n-m-1$. Let $M$ be the mapping cylinder of $\hl_s(X,Y)\to Y$, let $L^T$ be $L\cap Y$, and let $N^T$ be an NSDRN of $L^T$ in $L$.
Then 
$$IH^c_*(M, M-p^{-1}(y))\cong IH^c_{*-m}(cL, cN^T-y).$$
If $\bar p$ is a traditional perversity and the coefficient system is constant, then this is isomorphic to
\begin{equation*}\begin{cases}
IH^c_{*-m-1}(N^T), & *-m> k-\bar p(k+1),\\
\ker(IH^c_{*-m-1}(N^T)\to IH^c_{*-m-1}(L)), & *-m=k-\bar p(k+1),\\
IH^c_{*-m}(L,N^T),& *-m<k-\bar p(k+1).
\end{cases}
\end{equation*}
Also,
\begin{equation*}
IH^c_*(M-Y,M-(Y\cup p^{-1}(y))\cong IH^c_{*-m}(L-L^T, N^T-L^T).
\end{equation*}
\end{theorem}

\begin{corollary}\label{C: pm top}
In the preceding theorem, if $y$ lies in the top stratum of $Y$ (or any stratum that does not lie in the closure of another stratum), then $L^T$ is empty and so
$$IH^c_*(M, M-p^{-1}(y))\cong  IH^c_{*-m}(cL),$$ which is
\begin{equation*}
\begin{cases}
0, & *-m\geq k-\bar p(k+1),\\
IH^c_{*-m}(L),& *-m<k-\bar p(k+1),
\end{cases}
\end{equation*}
if $\bar p$ is a traditional perversity and the coefficient system is constant,
and
\begin{equation*}
IH^c_*(M-Y,M-(Y\cup p^{-1}(y)))\cong IH^c_{*-m}(L).
\end{equation*}
\end{corollary}

The proof of Theorem \ref{T: pm stalk} will be carried out over several lemmas and corollaries.

\begin{lemma}
Let $U$ be a distinguished neighborhood of $y$, and let $U^T=U\cap Y$. Then $$IH^c_*(M, M-p^{-1}(y))\cong IH^c_*(p^{-1}(U^T), p^{-1}(U^T)-p^{-1}(y))$$ and  
$$IH^c_*(M-Y,M-(Y\cup p^{-1}(y)))\cong IH^c_*(p^{-1}(U^T)-U^T, p^{-1}(U^T)-(p^{-1}(y)\cup U^T)).$$
\end{lemma}
\begin{proof}
This is just excision \cite[Lemma 2.11]{GBF10}.
\end{proof}

\begin{lemma}
Let $\mc M$ be the mapping cylinder of the restricted holink evaluation $\hl_s(U,U^T)\to U^T$. Then $(p^{-1}(U^T), p^{-1}(U^T)-p^{-1}(y))$ is stratum-preserving fiber homotopy equivalent to $(\mc M,\mc M-\td p^{-1}(y))$, where $\td p$ is the restriction of $p$ to $\mc M$. 
\end{lemma}
\begin{proof}
There is  an inclusion of the latter pair into the former. Let $\delta: U^T\to (0,\infty)$ be such that $B_{\delta}(x)\subset U$ for all $x\in U^T$. Then the mapping cylinder of $\hl^{\delta}_s(U,U^T)\to U^T$ is embedded in both space pairs, and it is stratum-preserving homotopy equivalent to each by the path shrinking arguments of \cite{Q1}.
\end{proof}

\begin{corollary}
$$IH^c_*(p^{-1}(U^T), p^{-1}(U^T)-p^{-1}(y))\cong IH^c_*(\mc M, \mc M-\td p^{-1}(y)),$$ and $$IH^c_*(p^{-1}(U^T)-U^T, p^{-1}(U^T)-(p^{-1}(y)\cup U^T))\cong IH^c_*(\mc M-U^T, \mc M-(U^T\cup \td p^{-1}(y))).$$ 
\end{corollary}
\begin{proof}
This follows from the preceding lemma, noting that all homotopy equivalences hold $U^T$ fixed and take the complement of $U^T$ to the complement of $U^T$.
\end{proof}

We can use the homeomorphism $U\cong \R^m\times cL$ to treat $L$ and $cL$ as specific subspaces of $U$ and hence of $X$. In particular, we choose $L$ so that the cone point of $cL$ is $y$. By the compatibility of the distinguished neighborhood structure with the stratification of $X$, $U^T$ has the form  $\R^m\times c(L\cap Y)$. Let $L^T$ denote $L\cap Y$ ($T$ for tangential). 

\begin{lemma}
$\mc M$ is stratum and fiber-preserving homotopy equivalent to $\R^m\times \td p^{-1}(cL^T)$. 
\end{lemma}
\begin{proof}
$\td p:\mc M\to U^T$ is a stratified fibration over $U^T\cong \R^m\times cL^T$. Since this base space is a product with the unstratified $\R^m$, the lemma follows as for ordinary fibration theory. The specific stratified analogue is given in \cite[Lemma 3.12]{GBF3}.
\end{proof}

\begin{corollary}
$$IH^c_*(\mc M, \mc M-\td p^{-1}(y))\cong IH^c_{*-m}(\td p^{-1}(cL^T), \td p^{-1}(cL^T-y)),$$ and $$IH^c_*(\mc M-U^T, \mc M-(U^T\cup \td p^{-1}(y)))\cong IH^c_{*-m}(\td p^{-1}(cL^T)-U^T, \td p^{-1}(cL^T-y)-U^T).$$ 
\end{corollary}
\begin{proof}
$(\mc M, \mc M-\td p^{-1}(y))\sim_{spfhe}(\R^m\times \td p^{-1}(cL^T), \R^m\times \td p^{-1}(cL^t)-\td p^{-1}(y))\sim_{spfhe} (R^m, R^m-\vec 0)\times (\td p^{-1}(cL^T), \td p^{-1}(cL^T-y)$. Thus the results follow from the K\"unneth theorem for intersection homology \cite{Ki}. Alternatively, one could proceed as in the proof of \cite[Proposition 2.20]{GBF10}.
\end{proof}

\begin{lemma}
$\td p^{-1}(cL^T)$ is stratum-preserving fiber homotopy equivalent to the mapping cylinder of $\hl_s(cL, c L^T)$. 
\end{lemma}
\begin{proof}
This will follow if we show that the restriction of $\pi:\hl_s(U,U^T)\to U^T$ to $\pi^{-1}(cL^T)\to cL^T$ is  stratum-preserving fiber homotopy equivalent to $\hl_s(cL,c L^T)\to cL^T$. Clearly the latter embeds in the former. But  $\pi^{-1}(cL^T)$ consists of paths in $U\cong \R^m\times cL$ that are stratum-preserving on $(0,1]$ and whose endpoints lie in $cL^T$. But there is clearly a stratum-preserving deformation retraction $\R^m\times cL\to cL$, and this retraction induces a deformation retraction from $\pi^{-1}(cL^T)$ to $\hl_s(cL,cL^T)$. Extending to the mapping cylinder completes the proof.
\end{proof}

\begin{corollary}
Let $\mf M$ be the mapping cylinder of $\hl_s(cL,cL^T)$ with projection $\check p:\mf M\to cL^T$. Then $$IH^c_{*-m}(\td p^{-1}(cL^T), \td p^{-1}(cL^T-y))\cong IH^c_{*-m}(\mf M,\check p^{-1}(cL^T-y))$$ and $$IH^c_{*-m}(\td p^{-1}(cL^T)-U^T, \td p^{-1}(cL^T-y)-U^T)\cong IH^c_{*-m}(\mf M- cL^T,\check p^{-1}(cL^T-y)- (cL^T-y)).$$ 
\end{corollary}

\begin{lemma}
$cL$ is an NSDRN of $cL^T$. If $N^T$ is an NSDRN of $L^T$ in $L$, then $cN^T-y$ is an NSDRN of $cL^T-y$.
\end{lemma}
\begin{proof}
The second statement is easy since $N^T$ is an NSDRN of $L^T$ by definition and $(cN^T-y,cL^T-y)\cong (N^T\times \R, L^T\times \R)$.

For the first part, continue to let $N^T$ be an NSDRN of $L^T$. We will write a nearly stratum-preserving deformation retraction of $cL$ to $L$ explicitly. The idea is to let points in $cN^T$ use their retraction paths to get to $cL^T$, while points further away simply retract to the cone point. Of course the issue is to interpolate properly between the two cases. For this,  let $\rho:L\to [0,1]$ be a continuous function such that $\rho(x)=1$ if and only if $x\in L^T$ and $\rho=0$ outside a neighborhood $U$ of $L^T$ such that $\bar U\subset N^T$. If $x\in N^T$, let $\omega_x$ be the retraction path of $x$ under a fixed nearly-stratum preserving deformation retraction of $N^T$ to $L^T$; then $\omega_x(1)=x$ and $\omega_x(0)\in L^T$. If $x\notin N^T$, we still let $\omega_x(1)=x$. We represent a point in $cL$ by the coordinates $(x,s)$, where $x\in L$ and $s\in [0,1]$. Of course $(x,0)=(y,0)$ for all $x,y\in L$. 

Then define
\begin{equation*}
H((x,s),t)=
\begin{cases}
(\omega_x(1-t), s(1-2t +2t\rho(x))),     & \rho(x)\geq 1/2, \\
(\omega_x(1-2t\rho(x)), s(1-t)),    & \rho(x)\leq 1/2 .
\end{cases}
\end{equation*}
This is well-defined, since for $x\notin N^T$, $\rho(x)=0$ and we only use $\omega_x(1)=x$ in the second expression. Each expression is continuous and  when $\rho(x)=1/2$ the expressions agree, so $H$ is continuous. For $t=0$,  $H((x,s),0)=(x,s)$, while at time $t=1$, the points from the first part land in $L^T$ and those from the second part land in the cone point. If $(x,s)\in (cL^T-y)$, then $\rho(x)=1$ and $H((x,s),t)=(\omega_x(1-t), s)=(x,s)$. The cone point stays fixed throughout. It is easy to check that this retraction is stratum-preserving.
\end{proof}

\begin{corollary}
$IH^c_{*-m}(\mf M,\check p^{-1}(cL^T-y))\cong IH^c_{*-m}(cL, cN^T-y)$, and $IH^c_{*-m}(\mf M- cL^T,\check p^{-1}(cL^T-y)- (cL^T-y))\cong IH^c_{*-m}(cL-cL^T, cN^T-cL^T)$. 
\end{corollary}
\begin{proof}
Let $\delta:cL^T-y\to (0,\infty)$ be such that $\delta(z)<\frac{1}{2} d(z, cL-cN^T)$. Let $A^{\delta}$ be the intersection of $\check p^{-1}(cL^T-y)$ with the mapping cylinder of the projection $\hl_s^{\delta}(cL,cL^T-y)\to cL^T-y$.  Then the inclusion $A^{\delta}\into \check p^{-1}(cL^T-y)$ is a homotopy equivalence by the arguments of \cite{Q1}. And we have maps $$(\mf M,\check p^{-1}(cL^T-y))\hookleftarrow(\mf M,A^{\delta})\to (cL, cN^T-y).$$ 
The left arrow is given by inclusion, and the right one is given by the modified map to retraction paths as in the stratum-preserving homotopy equivalence of \cite[Proposition A.1]{GBF3}. $\mf M\to cL$ is a stratum-preserving homotopy equivalence by \cite[Proposition A.1]{GBF3}, since $cL$ is an NSDRN of $cL^T$ by the preceding lemma. Similarly $A^{\delta}\to cN^T-y$ is a stratum-preserving homotopy equivalence, since it filters as 
$$A^{\delta}\into \text{mapping cylinder of }(\hl_s(cN^T-y,cL^T-y)\to cL^T-y)\to cN^T-y,$$ and each of these maps is a stratum-preserving homotopy equivalence. Making a diagram of the long exact sequences of pairs in intersection homology, we obtain   $IH^c_{*-m}(\mf M,\check p^{-1}(cL^T-y))\cong IH^c_{*-m}(cL, cN^T-y)$ by two applications of the five-lemma. 

The second statement follows the same way by removing the intersections with $Y$.
\end{proof}

\begin{lemma}
$IH^c_{*-m}(cL-cL^T, cN^T-cL^T)\cong IH^c_{*-m}(L-L^T, N^T-L^T)$, and if $\bar p$ is a traditional perversity and the coefficient system is constant,
\begin{equation*}
IH^c_{*-m}(cL, cN^T-y)=
\begin{cases}
IH^c_{*-m-1}(N^T), & *-m> k-\bar p(k+1),\\
\ker(IH^c_{*-m-1}(N^T)\to IH^c_{*-m-1}(L)), & *-m=k-\bar p(k+1),\\
IH^c_{*-m}(L,N^T),& *-m<k-\bar p(k+1).
\end{cases}
\end{equation*}
\end{lemma}
\begin{proof}
This follows from the cone formula, stratum-preserving homotopy equivalence, and some minor diagram chasing. 
\end{proof}

Putting together this string of lemmas and corollaries proves Theorem \ref{T: pm stalk}. \hfill \qedsymbol

\section{Indications of applications}\label{S: applications}

In this section, we mention briefly some elementary applications of the preceding spectral sequence machinery. 

The original motivation for this work was to be able to say something about the torsion of the intersection homology modules of regular neighborhoods, particularly neighborhoods of the embedding singularities of non-locally-flat knots \cite{GBF1}. For example, let  $X$ be a pseudomanifold and $N$ a regular neighborhood of a pure subset $Y$ of $X$. Suppose further that the intersection homology modules of the links of the points of $Y$ are all $\rho$-torsion, where $\rho$ is some element of the ground field $R$. Then it follows by applying standard arguments to our spectral sequences that the intersection homology modules of $N$ must also be $\rho$-torsion. Such arguments are used in \cite{GBF3, GBF5} to study the intersection Alexander invariants \cite{GBF2} of non-locally-flat PL knots. Although our results  allow for more general geometric situations than those of \cite{GBF3, GBF5}, we will not pursue this issue further here; we refer the reader to \cite{GBF3,GBF5} to see the patterns of such applications. Similar arguments are made in \cite{LM06} in Maxim's study of Alexander invariants of algebraic hypersurfaces.

As another sample application, let us consider the intersection homology of approximate tubular neighborhoods of locally-flat submanifolds in the topological category. So let $X$ be a topological $n$-manifold, and let $Y$ be an embedded locally-flat $m$-submanifold. The local-flatness condition ensures that the filtered space $X\supset Y$ is a topological pseudomanifold. Suppose that $Y$ has an approximate tubular neighborhoods $N$; it is shown in \cite{HTWW} that this will always be the case if $\dim X\geq 5$ (note that this paper predates some of the language of \cite{Hu02}, including the term ``approximate tubular neighborhood''). In this case, by Corollaries  \ref{C: NSDRN IHC}, \ref{C: manifold base PD}, and \ref{C: pm top}, the $E^2$ terms of the spectral sequence for $IH^c_*(N)$ are $H^c_{m-p}(Y; \mc{IH}^c_{n-m-q}(cL))$. But in this case, $cL\cong cS^{n-m-1}$. So if $\bar p(n-m)\geq 0$, then, using the cone formula, we see that the only nontrivial terms are $H^c_{m-p}(Y; \mc{IH}^c_{0}(cL))\cong  H^c_{m-p}(Y; \Z)$. So, the spectral sequence collapses, and, after properly organizing the indices, we see that $IH^c_*(N)\cong H^c_*(Y)$. Of course this is what one expects in this case, since the pair $X\supset Y$ constitutes a pseudomanifold, and intersection homology with traditional perversities is a topological invariant in this case. So $IH^c_*(N)\cong H_*(N)\cong H_*(Y)$. Note, however, that the spectral sequence computation at no point uses this topological invariance property.

On the other hand, if $\bar p(n-m)<0$ (in which case $\bar p$ must be a loose perversity \cite{Ki}), then the nontrivial $E^2$ terms are $H^c_{m-p}(Y; \mc{IH}^c_{0}(cL))\cong H^c_{m-p}(Y; \Z)$  and $H^c_{m-p}(Y; \mc{IH}^c_{n-m-1}(cL))\cong H^c_{m-p}(Y; \td \Z)$, where $\td Z$ is a possibly twisted coefficient system with $\Z$ stalks. In this case, the spectral sequence gives us the homology of the normal spherical fibration to $Y$ in $X$, which is also $H_*(N-Y)$. This is not completely obvious, however, because singular chains in $IC^c_i(X)$ with $i\geq n-m-\bar p(n-m)$ are allowed to intersect $Y$. 

We can also look at $IH^{\infty}_*(N)$. Putting together Theorem \ref{T: closed SS}, Propositions \ref{P: rel stalk} and \ref{P: easy pm}, and the arguments of Corollary \ref{C: manifold base PD}, there is a spectral sequence for $IH^{\infty}_{n-*}(N)$ with $E^2_{p,q}\cong H^{\infty}_{m-p}(Y; \mc{IH}^c_{n-m-q}(D^{n-m},S^{n-m-1}))$. If $p$ is a traditional perversity, the only nontrivial terms are $H^{\infty}_{m-p}(Y; \mc{IH}^c_{n-m}(D^{n-m},S^{n-m-1}))\cong H^{\infty}_{m-p}(Y; \td \Z)$. So if the associated normal spherical fibration of $Y$ in $X$ is untwisted then we get $IH^{\infty}_*(N)\cong H_{*-{n-m}}^{\infty}(Y)$, which we would certainly expect if $N\cong Y\times \R^{n-m}$. On the other hand, if $\bar p(n-m)<0$, then $IH^{\infty}_*(N)\cong 0$!

In forthcoming work, we will apply the spectral sequences developed here to study the local intersection homology groups on manifold homotopically stratified spaces. This will allow us to relate singular intersection homology to the Deligne sheaf construction \cite{GM2, Bo} on such spaces.

\section{Properties of approximate tubular neighborhoods}\label{S: ATN}

In this section, we provide the proofs of the properties of approximate tubular neighborhoods used in Section \ref{S: neighborhoods}.

We recall from \cite{Hu02} that $N$ is  an approximate tubular neighborhood of a pure subset $Y$ in $X$ if  there is a manifold stratified approximate fibration (MSAF) $p:N-Y\to Y\times \R$ such that the teardrop $(N-Y)\cup_pY$ is homeomorphic to $N$. 

An MSAF $q: A\to B$ is a proper map between manifold stratified spaces such that the following lifting condition is satisfied:
Given a diagram
\begin{diagram}
Z&\rTo^f & A\\
\dTo^{\times 0}&&\dTo_q\\
Z\times I&\rTo^F&B,
\end{diagram}
such that $Z$ is arbitrary and $F$ is a stratum-preserving homotopy, there is a \emph{weak stratified controlled solution} $\td F: Z\times I\times [0,1)\to A$ that is stratum-preserving along $I\times [0,1)$, satisfies $\td F(z,0,t)=f(z)$, and is such that the function $\bar F:Z\times I\times I\to B$ defined by $\bar F|Z\times I\times [0,1)=p\td F$ and $\bar F|Z\times I\times \{1\}=F$ is continuous. 

By \cite[p. 873]{Hu02}, if $N$ is an approximate tubular neighborhood, then the natural extension $\td p: N\to Y\times (-\infty,\infty]$ is also proper.

\begin{lemma}[Lemma \ref{L: ATN=NSDRN}]
Suppose that $N$ is an approximate tubular neighborhood of $Y$ in the MHSS $X$. Then there is a nearly stratum-preserving deformation retraction taking $N$ into $Y$ rel $Y$.
\end{lemma}
\begin{proof}
The proof of the lemma follows roughly the proof of Proposition 3.2 of \cite{Hug99}, which is concerned with a similar deformation retraction of the \emph{damped mapping cylinder} of $p'$, which is the composition of $p$ with the inclusion $Y\times \R \into Y\times (-\infty,\infty]$. We deviate from Hughes's proof only in that we restrict attention to $N$ rather than this mapping cylinder, which leads to some minor simplifications, and we fill in some details regarding the construction of the maps $\td g$. The reader is encouraged to compare with \cite{Hug99}.

For each $i=0,1,\ldots$, let $g^i:Y\times (-\infty,\infty]\times [\frac{i}{i+1},\frac{i+1}{i+2}]\to Y\times (-\infty,\infty]$ be an isotopy such that 
\begin{enumerate}
\item $g^i_{i/(i+1)}=\text{id}$,
\item $g^i$ is constant off of  $Y\times (-\infty, i+0.75]$,
\item $\text{Im}(g^i_{(i+1)/(i+2)})\subset Y\times [i+0.25,\infty]$,
\item $g^i$ is fiber preserving over $Y$.
\end{enumerate}

We will construct maps (analogous to Hughes's) $\td g^i:N\times [\frac{i}{i+1},\frac{i+1}{i+2}]\to N$, $i=0,1,\ldots$ such that 
\begin{enumerate}
\item\label{I: start} $\td g^i_{i/(i+1)}=\text{id}$,

\item\label{I: sp} $\td g^i$ is stratum-preserving along $[\frac{i}{i+1},\frac{i+1}{i+2}]$,

\item\label{I: iclose} $p\td g^i_t$ is $1/2^i$-close to $g^i_tp$ for each $t\in [i/(i+1), (i+1)/(i+2)]$,

\item \label{I: nclose} $p\td g_t^i|p^{-1}(Y\times [n,\infty))$ is $1/2^n$-close  to $g_t^ip|p^{-1}_{Y\times \R}(Y\times [n,\infty))$ for each $t\in[i/(i+1), (i+1)/(i+2)]$ and $n=0,1,\ldots$

\item \label{I: Y} For $y\in Y$, $g_t^i(y)=y$ for all $t\in [\frac{i}{i+1},\frac{i+1}{i+2}]$. 
\end{enumerate}

For this, let $\phi^i:N\times [\frac{i}{i+1},\frac{i+1}{i+2}]\times [0,1)\to N$ be a stratified controlled lift of $g^i|Y\times \R$ extending the identity map on $N-Y$. We define a function $\alpha: N\times [\frac{i}{i+1},\frac{i+1}{i+2}]\to (0,1]$ as follows: Let $d$ be the distance on $Y\times \R$. Let $\alpha$ be continuous and such that $d(p\phi^i(z,t,\alpha(s)),g^i_tp(z))< 1/2^i$ for all $z\in N-Y$ and $t\in [\frac{i}{i+1},\frac{i+1}{i+2}]$, and such that  $d(p\phi^i(z,t,\alpha(s)),g^i_tp(z))< 1/2^n$, for all $z\in p^{-1}(Y\times [n, \infty))$. It is possible to find such an $\alpha$ since we know that the function defined by $p\phi^i$ on  $N-Y\times [\frac{i}{i+1},\frac{i+1}{i+2}]\times [0,1)\to N-Y$ and $g^ip$ on $N-Y\times [\frac{i}{i+1},\frac{i+1}{i+2}]\times 1\to N-Y$
is continuous, by the definition of the lifting property of MSAFs. We also recall that $N$ is metric, hence paracompact. 

Define $\td g^i$ by 
\begin{equation*}
\td g^i(z,t)=
\begin{cases}
z,& z\in Y,\\
\phi^i(z,t,\alpha(z,t)),& z\in N-Y.
\end{cases}
\end{equation*}
It is clear by construction that properties \eqref{I: iclose}, \eqref{I: nclose}, and \eqref{I: Y} are satisfied. Properties \eqref{I: start} and \eqref{I: sp} are immediate from the definition of MSAFs, as is the continuity of $\td g^i$ on $N-Y\times [\frac{i}{i+1},\frac{i+1}{i+2}]$. Lastly, to see that $\td g^i$ is continuous on $Y$, we note that by the teardrop topology, points $y\in Y$ have arbitrarily small neighborhoods of the form $U\cup(p^{-1}(U)\times(M,\infty))$ for $U$ a neighborhood of $y$ in $Y$ and $M\in \R$. So if we take the arbitrary neighborhood $W=U\cup( p^{-1}(U)\times(M,\infty))$ for $U=B_{\epsilon}(y)$, where $B_{\epsilon}(y)$ is the ball of radius $\epsilon$ about $y$ in $Y$ and $\epsilon$ and $M>0$ are arbitrary, then $\td g^i$ maps the neighborhood $B_{1/2^n}(y)\cup( p^{-1}(B_{1/2^n}(y))\times (n,\infty))$ into $W$ if $\frac{1}{2^{n-1}}<\epsilon$ and $n>M+1$.  

Now define $\td g:N\times [0,1)\to N$ by $\td g_t=\td g^i_t\circ \td g^{i-1}_{i/(i+1)}\circ \cdots \circ \td g^0_{1/2}$ for $\frac{i}{i+1}\leq t\leq \frac{i+1}{i+2}$. Then $\td g$ is stratum-preserving along $[0,1)$, and using the properties of $\td g^i$, the image $p\td g_i(N)$ is in $Y\times [i+0.25-1/2^i,\infty]$. In particular, as $t$ increases, $p\td g_t$ herds $N-Y$ towards the $\infty$ end of $Y\times (-\infty,\infty]$, and thus $\td g_t$ herds $N$ towards $Y$. 

Let $p_{\infty}:N-Y\to Y\times \infty$ denote the composition of $p$ with the retraction of $Y\times (-\infty,\infty]$ to $\infty$. Define $\hat g: N\times I\to N$ by
\begin{equation*}
\hat g (x,t)=
\begin{cases}
x & x\in Y\\
\td g_t(x) & x\in N-Y, t\in [0,1)\\
\lim_{t\to 1} p_{\infty}\circ \td g_t(x) & x\in N-Y, t=1
\end{cases}
\end{equation*}

It is clear that $\hat g$ is well defined and continuous on $(N-Y)\times [0,1)$. It follows easily from the properties of $\td g$ and $g$ that it is also continuous on $N\times [0,1)$. We need to check what happens at time $1$. 

The well-definedness of $\hat g_1$ also follows: By the choice of $g^i$ as fiber-preserving over $Y$, from the definition of $\td g$, and from the properties of $\td g^i$, for all $\epsilon>0$ there exists $t_0<1$  such that $d_Y(p_{\infty}\td g_t,p_{\infty}\td g_s)<\epsilon$ if $t_0\leq s,t<1$, where $d_Y$ is the metric on $Y$ induced from that on $N$. 

The only remaining thing to verify is the continuity of $\hat g$ at time $1$. So let $(x,1)\in N\times [0,1]$, and let $y=\hat g_1(x)$. Neighborhoods of $y$ in $N$ have the form $W=p^{-1}(U)\cup(U\cap Y)$, where $U$ is a neighborhood of $y=y\times \infty$ in $Y\times (-\infty,\infty]$. So to check continuity of $\hat g$ at $(x,1)$, we need only check that if we are given such a $W$, then $p\hat g$ takes all points sufficiently close to $(x,1)$ in $N\times [0,1]$ to $W$. But this once again follows from the properties of $\td g^i$ and epsilon arguments similar to those used in showing the continuity of $\td g^i$. 
\end{proof}

\begin{corollary}[Corollary \ref{C: ATN=NSDRN}]
Approximate tubular neighborhoods in MHSSs are NSDRNs. 
\end{corollary}

\begin{proposition}[Proposition \ref{P: ATN=OST}]
Approximate tubular neighborhoods of pure subsets of MHSSs are outwardly stratified tame (see Definition \ref{D: ost}).
\end{proposition}
\begin{proof}

Let $N$ be an approximate tubular neighborhood of a pure subset $Y$ in $X$.
Since $X$ is metric, separable, and locally compact, so are $N$ and $Y$, and thus $Y$ is also $\sigma$-compact. Let $K_1,K_2, \ldots$ be a sequence of compact subsets of $Y$ such that $K_i\subset \text{int}(K_{i+1})$ and $Y=\cup_i K_i$. Let $C_0=Y$ and $C_i= Y\cup p^{-1}(K_i\times (-i,\infty))$ for $i>0$. Let $Z$ be a metric space, and let $f: Z\to N$ be a proper map such that $f(Z)\subset N-Y$. We must construct a proper stratum-preserving homotopy $H:Z\times [0,\infty)\to N$ such that $H(Z\times [0,\infty))\subset N-Y$.

We claim that we can construct a stratum-preserving homotopy $H:Z\times [0,\infty)\to N$ such that 
\begin{enumerate}
\item $H(\cdot,0)=f$,
\item\label{I: outside} $H(Z\times [i,\infty))\subset N- C_i$ for $i=0,1,\ldots$, and 
\item\label{I: fix} For $t\in[0,i]$ and  $z\in f^{-1}(N-\text{int}(C_{i+2}))$, $H(z,t)=f(z)$.
\end{enumerate}
Let us first show that this suffices by showing that such an $H$ is proper. The condition on the image remaining in $N-Y$ follows automatically from property \eqref{I: outside}.

Let $A$ be any compact set of $N$. We must show that $H^{-1}(A)$ is compact. Since $A$ is compact, so are $p(A)\subset Y\times \R$ and the projections of $p(A)$ to $Y$ and to $\R$. Thus $p(A)\subset K_j\times (-N,N)\subset Y\times \R$ for some $j$ and $N$. But this implies that $A\subset C_m$ for some $m\in \N$. So by property \eqref{I: outside} of $H$, $H^{-1}(A)\subset Z\times [0, m]$. But by property \eqref{I: fix}, if $f(z)\subset N-C_{m+2}$, then $H(z\times [0,m])\not\subset C_m$. Thus $H^{-1}(A)\subset f^{-1}(C_{m+2})\times [0,m]$, and it suffices to see that $f^{-1}(C_{m+2})$ is compact. 

By Lemma \ref{L: proper is closed}, below, $f(Z)$ is closed since $f$ is proper.  Furthermore, by \cite[p. 873]{Hu02}, the natural extension $\td p: N\to Y\times (-\infty,\infty]$ is proper, so again by Lemma \ref{L: proper is closed},  $\td p(f(Z))$ is closed. Since $f(Z)$ does not intersect $Y$, $\td pf(Z)$ does not intersect $Y\times \{\infty\}\subset Y\times (-\infty,\infty]$, and thus $\td p(f(Z)\cap C_{m+2})=p(f(Z)\cap C_{m+2})\subset K_{m+2}\times [-(m+2), M]$, for some $M$. But then, since $p$ is proper, $f(Z)\cap C_{m+1}$ is compact, and so, since $f$ is proper, $f^{-1}(C_{m+2})$ is compact.

It remains to show that we can find an $H$ satisfying the given conditions. We will construct $H$ inductively in pieces $H:Z\times [i-1,i]$, $i=1,2, \ldots$. The induction starts with $H|Z\times 0=f$.

 Let $h^i: Y\times [i-1,i]\to Y\times \R$ be a homotopy that retracts $Y\times \R$ into $Y\times (-\infty,-i-\frac{1}{2})$ as follows: 
\begin{equation*}
h^i((y,s),t)=
\begin{cases}
(y,s), & s\leq -i-\frac{1}{2},\\
(y, (-i-\frac{1}{2})(t-(i-1))+(1-t+(i-1))s), &s\geq -i-\frac{1}{2}.
\end{cases}
\end{equation*}
Clearly $h^i$ is stratum-preserving. 

Consider the diagram:
\begin{diagram}
Z&\rTo^{H(\cdot, i-1)} & N\\
\dTo^{\times (i-1)}&&\dTo_p\\
Z\times [i-1,i]&\rTo^{h^i(p\times \text{id}_{[i-1,i]})(H(\cdot,i-1)\times \text{id}_{[i-1,i]})} &Y\times \R.
\end{diagram}
Since $p$ is an MSAF, there is a stratum-preserving controlled lift $\td h^i:Z\times [i-1,i]\times [0,1)\to N$ such that $\td h^i(z,i-1,t)=H(z,i-1)$ and such that the function $\bar h^i:Z\times [i-1,i]\times I\to Y\times \R$ defined by $\bar h^i|Z\times [i-1,i]\times [0,1)=p\td h^i$ and $\bar h^i|Z\times [i-1,i]\times \{1\}=h^i(p\times \text{id}_{[i-1,i]})(H(\cdot,i-1)\times \text{id}_{[i-1,i]})$ is continuous. 

As already observed, $f(Z)\cap C_{m}$ must be compact for any $m$ since $f$ is proper. Similarly, $H(\cdot, i-1)$ will be proper by induction, using the same arguments presented above to show that all of $H$ will be proper if the conditions on $H$ are satisfied.  Thus $H(Z,i-1)\cap C_{m}$ is compact for any $m$, as is $H(\cdot, i-1)^{-1}(C_{m})$. 

Now consider the image of $T=H(Z,i-1)\cap C_{i+2}$ under $h^i(\cdot, i)p$.  We have $h^i(p(T),i)\subset Y\times (-\infty,-i-\frac{1}{2})$ and so $p^{-1}(h^i(p(T),i))\subset N-C_{i}$. In other words,\linebreak $p^{-1}(\bar h^i(H(\cdot, {i-1})^{-1}(C_{i+2}),i,1))\subset N-C_{i}$. Thus by continuity, there is an open neighborhood $W_1$ of $H(\cdot, i-1)^{-1}(C_{i+2})\times i\times 1$ in $Z\times [i-1,i]\times I$ such that the image of this neighborhood under $\bar h^i$ lies outside of $\td p C_i$. Since $H(\cdot, i-1)^{-1}(C_{i+2})$ is compact, there is an open neighborhood  $U_1$  of $i\times 1$ in $[i-1,i]\times I$ such that $H(\cdot, i-1)^{-1}(C_{i+2})\times U_1\subset W_1$.

We will use $U_1$ in ensuring that $f(Z)$ gets out of $C_i$, but we must also ensure that nothing creeps back into $C_{i-1}$ in the middle of the process. For this, consider that, again by induction, $H(Z,i-1)\subset N-C_{i-1}$. By the definition of $h^i$ and $\bar h^i$, $\bar h^i(Z\times [i-1,i]\times 1)\cap K_{i-1}\times [-(i-1),\infty)=\emptyset$. So there is an open neighborhood of $W_2$ of $H(\cdot, i-1)^{-1}(C_{i+2})\times [i-1,i]\times 1$ in $Z\times [i-1,i]\times I$ such that $\bar h^i$ of this neighborhood lies outside of $\td p C_{i-1}$. Let $U_2$ be an open neighborhood of $[i-1,i]\times 1$ in $[i-1,i]\times I$ such that $H(\cdot, i-1)^{-1}(C_{i+2})\times U_2\subset W_2$. 

Finally, we note that $\bar h^i(H(\cdot,i-1)^{-1}(C_{i+2}-\text{int}(C_{i+1}))\times[i-1,i]\times 1)\cap \td p C_i=\emptyset$, so there is an open neighborhood $W_3$ of  $H(\cdot,i-1)^{-1}(C_{i+2}-\text{int}(C_{i+1}))\times [i-1,i]\times 1$ such that $\bar h^i$ of this neighborhood lies outside of $\td pC_i$. Let $U_3$ be an open neighborhood of $[i-1,i]\times 1$ in $[i-1,i]\times I$ such that $H(\cdot, i-1)^{-1}(C_{i+2}-\text{int}(C_{i+1}))\times U_3\subset W_3$. 

Now, let $\gamma:[0,1]\to U_2\cap U_3\cap ([i-1,i]\times [0,1))$ be a path such that $\gamma(0)\in U_2\cap U_3\cap(i-1\times [0,1))$ and $\gamma (1)\in U_1\cap U_2\cap U_3\cap(i\times [0,1))$. 
Let $\rho:N-Y\to [0,1]$ be a continuous function such that $\rho^{-1}(1)=C_{i+1}-Y$ and $\rho^{-1}(0)=(N-Y)-\text{int}(C_{i+2})$.  And for $t\in[i-1,i]$, define $$H(z,t)=\td h^i(z,\gamma(\rho(H(z,i-1))(t-(i-1)))).$$ 

Let us see that this does what we want. By definition of $\gamma$, the image of $\gamma$ is in $[i-1,i]\times [0,1)$, so $H$ is well-defined and stratum-preserving because $\td h^i$ is. For $t\in[i-1,i]$ and $z$ such that $H(z,i-1)\in N-\text{int}(C_{i+2})$, $H(z,t)=\td h^i(z,\gamma(0))=\td h^i(z,i-1,I)=H(z,i-1)$, which is $f(z)$ by induction hypothesis.  

In particular, for $z$ such that $H(z,i-1)=f(z)\in N-\text{int}(C_{i+2})$, $H(z\times [i-1,i])=f(z)\subset N-C_i$. But also for $z\in Z$ such that $H(z,i-1)\in C_{i+2}$, $H(z,[i-1,i])\subset N-C_{i-1}$  because $\gamma$ lies in $U_2$, and $U_2$ was defined for this purpose. Hence $H(Z\times [i-1,i])\subset N-C_{i-1})$. 

For any point $z$ in $H(\cdot,i-1)^{-1}(C_{i+1})$, $\rho(H(z,i-1))=1$, so that $H(z,i)=\td h^i(z,\gamma(1))$, which lies outside of $C_i$ by the choice of $U_1$, and for any point $z$ in $H(\cdot,i-1)^{-1}(C_{i+2}-\text{int}(C_{i+1}))$, $\td h^i(z,\gamma(\rho(H(z,i-1))t))$ is outside of $C_i$ for all $t\in[0,1]$ by the choice of $U_3$. Thus $H(Z,i)\subset N-C_{i}$.
\end{proof} 

The following lemma is not difficult and may be well known, but since the author could not find it in the most standard texts on general topology, we include a proof here.

\begin{lemma}\label{L: proper is closed}
A proper map $f: X\to Y$ from a metric space to a locally-compact metric space is a closed map. 
\end{lemma}
\begin{proof}
Suppose that $f$ is not closed. Then there is a closed set $Z\subset X$ such that $f(Z)$ is not closed, and there exists a point $y\in Y$ such that $y\notin f(Z)$ but $y\in \overline{f(Z)}$. Thus there is a sequence of points $y_i \in f(Z)$ such $\lim y_i=y$. Let $K$ be a compact neighborhood of $y$ in $Y$. Without loss of generality, we can assume that $y_i\in K$ for all $i$. Choose points $z_i\in Z$ such that $f(z_i)=y_i$. Since $f$ is proper, $f^{-1}(K)\cap Z$ is compact and the sequence $z_i$ has a convergent subsequence $z_{i_j}$ converging to a point $z\in Z$. By continuity of $f$, $f(z)=y$, a contradiction.    
\end{proof}

\bibliographystyle{amsplain}
\bibliography{bib}

\end{document}